\documentclass[journal]{IEEEtran}

\usepackage{xspace,amssymb,epsfig,subfigure,syntonly,bm,bbm}
\usepackage{epsfig,amsmath,color}
\usepackage{xspace,syntonly}

\usepackage{algorithm,algorithmic}

\newcommand{\utwi}[1]{\mbox{\boldmath $#1$}}

\newcommand{\diag}{{\textrm{diag}}}

\newcommand{\cD}{{\cal D}}

\newcommand{\cL}{{\cal{L}}}
\newcommand{\cN}{{\cal N}}
\newcommand{\cP}{{\cal P}}

\newcommand{\cS}{{\cal S}}

\newcommand{\cT}{{\cal T}}
\newcommand{\cC}{{\cal C}}
\newcommand{\cE}{{\cal E}}
\newcommand{\cI}{{\cal I}}
\newcommand{\cF}{{\cal F}}
\newcommand{\cB}{{\cal B}}

\newcommand{\cH}{{\cal H}}

\newcommand{\cW}{{\cal W}}

\newcommand{\ba}{{\bf a}}

\newcommand{\be}{{\bf e}}
\newcommand{\bg}{{\bf g}}

\newcommand{\bs}{{\bf s}}

\newcommand{\bv}{{\bf v}}
\newcommand{\bi}{{\bf i}}

\newcommand{\bA}{{\bf A}}

\newcommand{\bI}{{\bf I}}

\newcommand{\bZ}{{\bf Z}}

\newcommand{\bY}{{\bf Y}}
\newcommand{\bV}{{\bf V}}

\newcommand{\bnu}{{\utwi{\nu}}}

\newcommand{\bPsi}{{\utwi{\Psi}}}

\newcommand{\bPhi}{{\utwi{\Phi}}}
\newcommand{\bxi}{{\utwi{\xi}}}
\newcommand{\brho}{{\utwi{\rho}}}

\newcommand{\bchi}{{\utwi{\chi}}}
\newcommand{\bmu}{{\utwi{\mu}}}

\newcommand{\biota}{{\utwi{\iota}}}
\newcommand{\bsigma}{{\utwi{\sigma}}}

\newcommand{\bbmA}{\bar{\bm{A}}}
\newcommand{\bbmM}{\bar{\bm{M}}}
\newcommand{\bbmZ}{\bar{\bm{Z}}}
\newcommand{\tbmZ}{\tilde{\bm{Z}}}
\newcommand{\bbmQ}{\bm{Q}}

\begin{document}

\newtheorem{definition}{Definition}
\newtheorem{remark}{Remark}
\newtheorem{proposition}{Proposition}
\newtheorem{lemma}{Lemma}


\title{Sparsity-leveraging Reconfiguration \\ of Smart Distribution Systems}

\author{\emph{Emiliano Dall'Anese and Georgios B. Giannakis}
\thanks{\protect\rule{0pt}{0.5cm}%
Submitted to IEEE Transactions on Power Delivery on March 31, 2013; revised on August 17, 2013. 
This work was supported by the Inst. of Renewable Energy and the
Environment (IREE) grant no. RL-0010-13, Univ. of Minnesota. Authors are with the Digital Technology Center and the Dept. of ECE, University of Minnesota, 200 Union Street SE, Minneapolis, MN 55455, USA. 
E-mails: {\tt \{emiliano, georgios\}@umn.edu} 
%
%
%
}
}

\markboth{}%
{Dall'Anese and Giannakis: Sparsity-leveraging Reconfiguration of Smart Distribution Systems }

\maketitle

\maketitle

\begin{abstract}

A system reconfiguration problem is considered for three-phase power distribution networks featuring distributed generation. In lieu of binary line selection variables, the notion of group sparsity is advocated to re-formulate the \emph{nonconvex} distribution system reconfiguration (DSR) problem into a \emph{convex} one. 
Using the duality theory, it is shown that the line selection task boils down to a shrinkage and thresholding operation on the line currents. Further, numerical tests illustrate the ability of the proposed scheme to identify meshed, weakly-meshed, or even radial configurations by adjusting a sparsity-tuning parameter in the DSR cost. Constraints on the voltages are investigated, and incorporated in the novel DSR  problem to effect voltage regulation. 
\end{abstract}


\section{Introduction}
\label{sec:Introduction}

The fundamental objective of distribution system reconfiguration (DSR) schemes is to identify the topology of a distribution network that is optimal in a well defined sense~\cite{Merlin75,Baran89,Shirmohammadi89}. DSR byproducts include balancing the network load~\cite{Chin02}, increasing the system security, and prompt (possibly network-wide) power delivery restoration in case of localized network failures. Computationally-affordable DSR schemes are increasingly advocated in modern distribution networks to enhance their efficiency and stability in the presence of distributed generation (DG), energy storage devices, as well as dispatchable and elastic loads~\cite{Hatziargyriou-PESMag}. 

Changes in the network topology are effected by opening or closing tie and sectionalizing line switches.
These switching operations can be either performed manually in situ, or, commanded remotely by a network controller. Thus, the DSR task is traditionally approached by associating a binary selection variable with each switch~\cite{Cho93,Schmidt95,Gomes06,Huang02,Khodr09,Moradzadeh12}. 
Unfortunately, this choice renders the resultant topology selection problem~\emph{NP-hard}~\cite{Nemhauser88}, and thus challenging to solve optimally and efficiently. This explains why heuristic schemes have been largely employed to select the status of the switches. For example, all switches are initially assumed closed in e.g.,~\cite{Merlin75,Gomes06}, and then some of them are progressively opened until a radial configuration is obtained. A search over relevant radial configurations based on approximate power flow methods is advocated in~\cite{Baran89}. Alternative methods rely on fuzzy multi-objective~\cite{Das06}, branch-and-bound~\cite{Cho93}, and genetic algorithms~\cite{Huang02}. Off-the-shelf solvers for mixed-integer linear programs were employed in conjunction with Bender's decomposition in~\cite{Khodr09}, to solve a joint DSR and optimal power flow (OPF) for balanced systems. Newton methods and branch-selection heuristic techniques were employed in~\cite{Schmidt95}. However, these schemes are tailored for balanced systems, and their computational complexity may become prohibitive as the size of the network increases. An efficient exhaustive search algorithm was proposed in~\cite{Morton00} to find the optimal radial configuration, based on successive tree transformations. However, this approach cannot be utilized when the sought topology is (weakly-)meshed~\cite{Cicoria04}. Finally, the heuristic of~\cite{Merlin75,Gomes06} was extended to the case of unbalanced systems in~\cite{Zidan11}. 

The present paper leverages contemporary compressive sampling tools~\cite{YuLi06,Wiesel11} to bypass binary optimization variables, and formulate a novel \emph{convex} DSR problem. Specifically, the formulated DSR problem is a second-order cone program (SOCP), and it can be solved efficiently even for distribution networks of large size and with densely deployed line switches, using primal-dual interior point methods tailored to SOCPs~\cite{Nesterov94,Lobo98}. DSR solvers able to find a new configuration in a few seconds (or even less that one second, as shown in the numerical test cases) are instrumental for network operators to quickly re-shape the distribution grid in case of localized outages~\cite{Pilo99}, and to gauge optimal topologies in case of abrupt load or (renewable-based) generation variations. \textcolor{black}{Different from DSR approaches applicable to balanced distribution networks~\cite{Baran89,Shirmohammadi89,Gomes06,Khodr09,Moradzadeh12,Das06,Morton00}, the formulation herein  accounts also for unbalanced loads and non-zero off-diagonal entries of line admittance matrices. }

The proposed convex formulation hinges on the notion of \emph{group-sparsity}, an underlying attribute of the currents flowing on the phases of distribution lines equipped with switches. This group-sparse problem structure allows one to discard binary optimization variables, and select the states of the switches by augmenting the DSR cost with a convex sparsity-promoting regularization function~\cite{YuLi06,Wiesel11}. As in conventional (group) sparse linear regression, it is shown here that the line selection task boils down to a shrinkage and thresholding operation. This is further corroborated through  experiments on a modified version IEEE 37-node feeder~\cite{testfeeder} and other two test systems, where a meshed, weakly-meshed, or radial configuration is obtained by simply adjusting a sparsity-tuning parameter. 

Unfortunately, PQ loads and DG units (modeled as PQ loads as well) involve challenging nonlinear power flow relations. However, since the aim here is to develop a DSR scheme that is computationally efficient and yet able to reliably discard lines involving high active power losses, the approximate load model developed in~\cite{Bolognani13} is advocated, and tailored to the three-phase setup.
Although this load model introduces an approximation error (which  becomes negligible for sufficiently large values the nominal voltage as shown analytically in~\cite{Bolognani13}, and further corroborated numerically here), the payoff here is huge, since  a \emph{convex} DSR problem can be formulated even in the presence of PQ loads.\footnote{{\it Notation:} Upper (lower) boldface
letters will be used for matrices (column vectors); $(\cdot)^\cT$ for transposition; $(\cdot)^*$ complex-conjugate; and, $(\cdot)^\cH$ complex-conjugate transposition; $\Re\{\cdot\}$ denotes the real part, and $\Im\{\cdot\}$ the imaginary part; $j = \sqrt{-1}$ represents the imaginary unit; and $\mathbb{I}_{\{\cdot\}}$ is the indicator function
($\mathbb{I}_{\{x\}}=1$ if $x$ is true, and zero otherwise). $|\cP|$ denotes the cardinality of set $\cP$; $\mathbb{R}^{N}$ and $\mathbb{C}^{N}$ represent the space of the $N \times 1$ real and complex vectors, respectively. Given a vector $\bv$ and a matrix $\bV$, $[\bv]_{\cP}$ denotes a $|\cP| \times 1$ sub-vector containing entries of $\bv$ indexed by the set $\cP$, and $[\bV]_{\cP_1,\cP_2}$ the $|\cP_1| \times |\cP_2|$ sub-matrix with row and column indexes described by $\cP_1$ and $\cP_2$. Further, $\|\bv \|_2 := \sqrt{\bv^{\cT}\bv}$ denotes the $\ell_2$ norm of $\bv$. Finally, $\mathbf{0}_{M\times N}$ and $\mathbf{1}_{M\times N}$ denote $M \times N$ matrices with all zeroes and ones, respectively.}

\section{Preliminaries and problem formulation}
\label{sec:Formulation}

Consider a portion of the power distribution grid located downstream of the distribution substation, that supplies
a number of industrial and residential loads, and may include DG. 
The considered three-phase network is modeled also as a directed graph\footnote{The symbols defined throughout the paper are recapitulated in Table~\ref{tab:nomenclature}.}  $(\cN, \cE)$, where $N$ nodes are collected in the set $\cN := \{1,\ldots,N\}$, and overhead or underground lines are represented by the set of (directed) edges $\cE := \{(m,n)\} \subset \cN  \times \cN$. 
Let node $1$ represent the point of common coupling (PCC), taken to be the distribution substation. Distribution systems typically have tie and sectionalizing switches, whose states (closed or open) determine the topology of the network. Thus, let $\cE_{R} \subset \cE$ collect the branches equipped with controllable switches. 

Let $\cP_{mn} \subseteq \{a_{mn},b_{mn},c_{mn}\}$ and $\cP_{n} \subseteq \{a_{n},b_{n},c_{n}\}$ denote the set of phases of line $(m,n) \in \cE$ and node $n \in \cN$, respectively; $I_{mn}^{\phi} \in \mathbb{C}$  the complex current flowing from node $m$ to node $n$ on phase $\phi$; $I_m^{\phi} \in \mathbb{C}$ the current injected at node $m \in \cN$ and phase $\phi \in \cP_m$; and, $V_m^{\phi} \in \mathbb{C}$ the complex line-to-ground voltage at the same node and phase. Lines $(m,n) \in \cE$ are modeled as $\pi$-equivalent components~\cite[Ch.~6]{Kerstingbook}, and the $|\cP_{mn}| \times |\cP_{mn}|$ phase impedance  matrix is denoted by $\bZ_{mn} \in \mathbb{C}^{|\cP_{mn}| \times |\cP_{mn}|}$.  Matrix $\bZ_{mn}$ is symmetric, full-rank, and it is obtained from the line primitive impedance matrix via Kron reduction~\cite[Ch.~4]{Kerstingbook}. Using $\bZ_{mn}$, the $|\cP_{mn}| \times 1$ vector $\bi_{mn} := [\{I_{mn}^{\phi}\}_{\phi \in \cP_{mn}}]^\cT$ collecting the currents flowing on each phase of line $(m,n) \in \cE$ can be expressed as 
\begin{align}
\bi_{mn} = \bZ_{mn}^{-1} \left( [\bv_m]_{\cP_{mn}} - [\bv_n]_{\cP_{mn}} \right) 
\label{line_currents}
\end{align}
with $\bv_{m} := [\{V_{m}^{\phi}\}_{\phi \in \cP_{m}}]^\cT$. 
Let $\bi_{n} := [\{I_{n}^{\phi}\}_{\phi \in \cP_{n}}]^\cT$ be the vector collecting the currents injected at node $n$, and $\{\be_n^{\phi}\}_{\phi \in \cP_{n}}$ and $\{\be_{mn}^{\phi}\}_{\phi \in \cP_{mn}}$ the canonical bases of $\mathbb{R}^{|\cP_n|}$ and $\mathbb{R}^{|\cP_{mn}|}$, respectively. Further, per-node $n \in \cN$, define the $|\cP_{n}| \times |\cP_{mn}|$ matrix $\bA_{mn}^{(n)} := \sum_{\phi \in \cP_{n}} \mathbb{I}_{\{\phi \in \cP_{mn}\}}  \be_{n}^{\phi} (\be_{mn}^{\phi})^{\cT}$. Suppose for brevity, that the entries of the line shunt admittance matrix are negligible (they are, in fact, typically on the order of $10-100$ micro Siemens per mile~\cite{testfeeder}). However, at the expense of minimally increasing complexity, perceptible effects of shunt admittances can be readily accounted for in the ensuing problem formulations. Under these conditions, and using the definition of $\bA_{jn}^{(n)}$, Kirchhoff's current law can be written per node $n$ as
\begin{align}
\bi_{n} + \sum_{j \in \cN_{\rightarrow n}} \hspace{-.2cm} \bA_{jn}^{(n)} \, \bi_{jn}  - \sum_{k \in \cN_{n \rightarrow}}  \hspace{-.2cm} \bA_{nk}^{(n)} \, \bi_{nk} = \mathbf{0}_{|\cP_n| \times 1}
\label{kcl}
\end{align}
where $\cN_{\rightarrow n} := \{j : (j,n) \in \cE\}$ and $\cN_{n \rightarrow} := \{k : (n,k) \in \cE\}$, respectively. 
Clearly, $\bi_{n} = \mathbf{0}_{|\cP_n| \times 1}$ if neither loads nor DG units are connected at node $n$.  

\begin{table}[t]
\caption{\textcolor{black}{Nomenclature and main definitions}}
\begin{center}
\textcolor{black}{
\begin{tabular}{ll}
 \hline  \\ 
$\cN$ & Set collecting the nodes of the distribution system \\
$\cE$ & Set collecting the distribution branches  \\
$\cE_R$ & Subset of branches equipped with controllable switches  \\
$x_{mn}$ & Variable indicating the status of switch $(m,n) \in \cE_R$ \\  
$\cP_n$ & Set of phases at node $n$, $\cP_n \subseteq \{a, b, c\}$  \\
$\cP_{mn}$ & Set of phases of branch $(m,n)$, $\cP_{mn} \subseteq \{a, b, c\}$  \\
$\cS$ & Subset of nodes featuring distributed generation  \\
$\bZ_{mn}$ & Phase impedance matrix of line $(m,n)$ \\
$\bY^{(s)}_{mn}$ & Shunt admittance matrix of line $(m,n)$ \\
$V_n^\phi$ & Complex line-to-ground voltage at phase $\phi$ of node $n$   \\
$I_n^\phi$ & Complex current injected at phase $\phi$ of node $n$   \\
$I_{mn}^\phi$ & Complex current on phase $\phi$ of line $(m,n)$    \\
$S_{mn}$ &  Complex power injected on line $(m,n)$ from node $m$ \\
$S_{L,n}^{\phi}$ & Complex load demanded at node $n$ on phase $\phi$  \\
$S_{G,n}^{\phi}$ & Complex power supplied at node $n$ on phase $\phi$  \\
$I_{mn}^{\textrm{max}}$ & Maximum value for $|I_{mn}^\phi|^2$  \\
$V_N$ & Nominal line-line voltage  \\ 
$\varphi_N^\phi$ & Nominal angle of phase $\phi \in \{a, b, c\}$ \\
$\bI_{k}$ & $k \times k$ identity matrix \\
$\{\be_{n}^{\phi}\}_{\phi \in \cP_{n}}$ & Canonical bases of $\mathbb{R}^{|\cP_n|}$ \\ 
$\{\be_{mn}^{\phi}\}_{\phi \in \cP_{mn}}$ & Canonical bases of $\mathbb{R}^{|\cP_{mn}|}$ \\ 
$\bi_{mn}$ & $\bi_{mn} := [\{I_{mn}^\phi\}_{\phi \in \cP_{mn}}]^\cT \in \mathbb{C}^{|\cP_{mn}|} $ \\
$\bxi_{mn}$ & $\bxi_{mn} := [\Re^\cT\{\bi_{mn}\}, \Im^\cT\{\bi_{mn}\}]^{\cT} \in  \mathbb{R}^{2 |\cP_{mn}|} $ \\
$\bi_{n}$ & $\bi_{n} := [\{I_{n}^\phi\}_{\phi \in \cP_{n}}]^\cT \in \mathbb{C}^{|\cP_{n}|} $ \\
$\biota_{mn}$ & $\bxi_{mn} := [\Re^\cT\{\bi_{n}\}, \Im^\cT\{\bi_{n}\}]^{\cT} \in  \mathbb{R}^{2 |\cP_{n}|} $ \\
$\cN_{\rightarrow n} $ &  $\cN_{\rightarrow n} := \{j : (j,n) \in \cE\}$ \\
$\cN_{n \rightarrow} $ & $\cN_{n \rightarrow} := \{k : (n,k) \in \cE\}$ \\
$\bv_{n}$ & $\bv_{n} := [\{V_{n}^{\phi}\}_{\phi \in \cP_{n}}]^\cT \in  \mathbb{C}^{|\cP_{n}|}$ \\ 
$\bA_{mn}^{(n)}$ & $\bA_{mn}^{(n)} := \sum_{\phi \in \cP_{n}} \mathbb{I}_{\{\phi \in \cP_{mn}\}}  \be_{n}^{\phi} (\be_{mn}^{\phi})^{\cT}$ \\
$\bbmZ_{mn} $ & $\bbmZ_{mn} := \bI_{2} \otimes \Re\{\bZ_{mn}\}$ \\
$\bbmA_{mn}^{(m)}$ & $\bbmA_{mn}^{(m)} := \bI_{2} \otimes \bA_{mn}^{(m)}$ \\
$\bbmM_{mn}^{\phi}$ & $\bbmM_{mn}^{\phi} := \bI_{2} \otimes \be_{mn}^\phi (\be_{mn}^\phi)^\cT$ \\
$\bPhi_n$ & $\bPhi_n := \diag(\{e^{\varphi_N^\phi}\}_{\phi \in \cP_n})$ \\
$\bsigma_{L,n}$ &  $\bsigma_{L,n} := [\Re^\cT\{[\{S_{L,n}\}]\}, \Im^\cT\{[\{S_{L,n}\}]\}]^\cT$ \\ 
$\bsigma_{G,n}$ &  $\bsigma_{G,n} := [\Re^\cT\{[\{S_{G,n}\}]\}, \Im^\cT\{[\{S_{G,n}\}]\}]^\cT$ \\ 
\hline 
\end{tabular}}
\end{center}
\label{tab:nomenclature}
\end{table}%

Let $S_{mn} :=  \bi_{mn}^{\cH} [\bv_m]_{\cP_{mn}}$ denote the total complex power injected on line $(m,n)$ from node $m$. If no power is dispelled through the distribution line $(m,n) \in \cE$, $S_{mn}$ coincides with the total power transferred to node $n$; that is, $S_{mn} = - S_{nm}$. A necessary condition for this identity to hold is to have an identically zero line impedance matrix. In fact, it readily follows from~\eqref{line_currents} that
$S_{mn} + S_{nm} = \bi_{mn}^{\cH}\bZ_{mn} \bi_{mn}$, and thus the total active power loss on line $(m,n) \in \cE$ amounts to 
\begin{align}
\hspace{-.2cm} \Delta P_{mn} := \Re\{S_{mn} + S_{nm}\}  & = \Re^\cT \{\bi_{mn}\} \Re\{\bZ_{mn}\} \Re\{\bi_{mn}\}  \nonumber \\
 & \hspace{-0cm} + \Im^\cT\{\bi_{mn}\} \Re\{\bZ_{mn}\} \Im\{\bi_{mn}\} . \label{Ploss} 
\end{align} 
However, since typical values of $\bZ_{mn}$ in overhead and underground distribution segments render $\Delta P_{mn}$ not negligible~\cite{testfeeder,Kerstingbook}, it is desirable to select the topology (meaning the states of switches on lines $\cE_{R}$), and adjust the complex line currents $\{\bi_{mn}\}$  traversing the selected lines, so that the overall real power loss $\sum_{(m,n) \in \cE} \Delta P_{mn}$ is minimized~\cite{Merlin75, Baran89}. 

Similar to various DSR renditions~\cite{Cho93,Schmidt95,Gomes06,Huang02,Khodr09,Moradzadeh12}, the topology selection problem will be first formulated in the ensuing subsection using binary line selection variables. However, since lines may be non-transposed and the spacings between conductors may be non-equilateral~\cite{Kerstingbook}, the off-diagonal elements of $\bZ_{mn}$ are non-zero~\cite{testfeeder,Kerstingbook}. Thus, different from DSR schemes tailored to balanced distribution networks (as in e.g.,~\cite{Baran89,Gomes06,Khodr09,Moradzadeh12}), the problem formulated here is able to capture current-coupling effects on the distribution lines. 

\subsection{Plain-vanilla DSR formulation}
\label{sec:dsr_integer}

Suppose for the moment that loads are modeled using ideal current generators (that absorb current from the network). PQ loads and DG power injections (which follow a constant PQ model as well~\cite{LoadModel}) will be considered in Section~\ref{sec:load_approximation}. Similar to~\cite{Cho93,Schmidt95,Gomes06,Huang02,Khodr09,Moradzadeh12}, associate a binary variable $x_{mn} \in \{0,1\}$ with line $(m,n) \in \cE_{R}$, and suppose that this distribution segment is used to deliver power to the loads if $x_{mn} = 1$. In this case, the DSR problem can be formulated as follows [cf.~\eqref{kcl}] 
\begin{subequations}
\label{recon_nonconvex_current}
\begin{align}
& \hspace{-1.8cm} (P1)    \min_{\{\bi_{mn}\}, \bi_{1}, \{x_{mn}\} }  \sum_{(m,n) \in \cE} \Re^\cT\{\bi_{mn}\} \Re\{\bZ_{mn}\} \Re\{\bi_{mn}\} \nonumber \\ 
&  \hspace{2.0cm} + \Im^\cT\{\bi_{mn}\} \Re\{\bZ_{mn}\} \Im\{\bi_{mn}\}   \label{P1cost}  \\
\textrm{subject~to~~}  & \eqref{kcl},~\mathrm{and} \nonumber \\
& \hspace{-1.3cm}  \Re^2\{[\bi_{mn}]_{\phi}\} + \Im^2\{[\bi_{mn}^{\Im}]_{\phi}\} \leq I_{mn}^{\textrm{max}}, \,\, (m,n) \in \cE \backslash \cE_{R}  \label{P1magnitude}  \\
& \hspace{-1.3cm}  \Re^2\{[\bi_{mn}]_{\phi}\} + \Im^2\{[\bi_{mn}^{\Im}]_{\phi}\} \leq I_{mn}^{\textrm{max}} x_{mn}, \, (m,n) \in \cE_{R}  \label{P1close}  \\
& \hspace{1.9cm} x_{mn}  \in \{0,1\} \, , \quad \,  (m,n) \in \cE_{R} \label{P1integer}
\end{align}
\end{subequations}
where constraints~\eqref{P1magnitude}--\eqref{P1close} are enforced on each phase $\phi \in \cP_{mn}$, and $I_{mn}^{\textrm{max}} \geq 0$ is a cap for $|I_{mn}^{\phi}|^2$. Clearly, when $x_{mn} = 0$, $\bi_{mn}$ is forced to zero by~\eqref{P1close},  thus implying that line $(m,n) \in \cE_R$ is not used (see also~\cite{Cho93,Schmidt95,Khodr09,Moradzadeh12}). When the desired topology is  radial, additional constraints are present~\cite{Khodr09, Moradzadeh12}. In particular, suppose that the graph $(\cN,\cE)$ contains $N_I$ cycles, and collect in the set $\cC_i$ the lines forming cycle $i = 1,\ldots, N_I$. Then, to obtain a tree network, it suffices to add in $(P1)$ the constraint $\sum_{(m,n) \in \cC_i} x_{mn} \leq |\cC_i| - 1 $ per cycle $i$. 

Matrices $\{\Re\{\bZ_{mn}\} \in \mathbb{R}^{|\cP_{mn}| \times |\cP_{mn}|}\}$ are typically positive definite and full-rank (see e.g., the real test cases in~\cite{testfeeder}); thus, it follows that the DSR cost~\eqref{P1cost} is strictly convex. However, presence of the binary variables $\{x_{mn}\}$ renders $(P1)$ a mixed-integer quadratic program (MIQP), which is \emph{nonconvex} and \emph{NP-hard}~\cite{Nemhauser88}. Finding its global optimum requires solving a number of subproblems (one per switch status) that increases exponentially ($2^{|\cE_{R}|}$) in the number of switches. This explains why heuristic schemes have been largely employed to select the network topology~\cite{Baran89,Cho93,Gomes06,Das06}. Alternatively, off-the-shelf solvers for mixed-integer (non)linear programs~\cite{Khodr09} and genetic algorithms~\cite{Huang02} have been employed.  
However, since these solvers are in general computationally-heavy, they are not suited for optimizing the operation of medium- and large-size distribution networks, or, for finding a post-outage system configuration in order to efficiently restore power delivery network-wide~\cite{Pilo99}.  

In the ensuing Section~\ref{sec:dsr_sparsity}, compressive sampling tools will be advocated to bypass binary selection variables, and re-formulate the DSR problem into a \emph{convex} one. But first, an approximate yet powerful load model is outlined next, and other possible cost functions are described in Section~\ref{sec:cost_functions}. 

\subsection{Approximate load model} 
\label{sec:load_approximation}

Consider the well-established ``exponential model'' relating injected (or supplied) powers with voltages $\{V_n^{\phi}\}$ and currents $\{I_n^{\phi}\}$ (see e.g.,~\cite{LoadModel})  
\begin{align}
\label{loadmodel}
V_n^{\phi} (I_n^{\phi})^* = S_n^\phi \left| \frac{\sqrt{3}}{V_N} V_n^{\phi} \right|^{\kappa_n} , \quad \phi \in \cP_n, n \in \cN \backslash \{1\}
\end{align}
where $V_N$ is the nominal line-line voltage magnitude of the system (e.g., $4.8$ kV for the IEEE 37-node feeder~\cite{testfeeder}); $S_n^\phi$ is the net complex power that would be injected on phase $\phi$ of node $n$ if $V_n^{\phi}$ were equal to the nominal voltage $\frac{V_N}{\sqrt{3}} e^{j \varphi^\phi_N}$, $\varphi^\phi_N \in \{0^\circ, -120^\circ, 120^\circ\}$; and, $\kappa_n \in \{0, 1, 2\}$ is a model parameter. Specifically, constant PQ, constant current, and constant impedance loads are obtained by setting $\kappa_n = 0$, $\kappa_n = 1$, and $\kappa_n = 2$, respectively. 

DG units are typically modeled as constant PQ loads (that supply power); thus, let $\cS \subset \cN$ be the (sub)set of nodes where DG units are present, and $\bs_{G,n} := [\{ S_{G,n}^{\phi }\}]^\cT$ the vector collecting the complex power $S_{G,n}^{\phi }$ supplied by these units on each phase of node $n$. Notice that multiple DG units may be present at each node; if this is the case, $S_{G,n}^{\phi }$ can be readily replaced by $\sum_{u = 1}^{N_U} S_{G,n,u}^{\phi }$, with $N_U$ the number of DG units at node $n$. Per phase $\phi \in \cP_{n}$, let $S_{L,n}^{\phi}$ denote the complex power demanded by a wye-connected load at the bus $n$. Finally, define $\bs_{L,n} := [\{ S_{L,n}^{\phi }\}]^\cT$, and suppose as usual that 
the voltages at the substation  $\bv_1 :=  [\frac{V_N}{\sqrt{3}} e^{j 0^\circ},\frac{V_N}{\sqrt{3}} e^{j -120^\circ}, \frac{V_N}{\sqrt{3}} e^{j 120^\circ}]^{\cT}$ are taken as reference for the phasorial representation~\cite{Kerstingbook}. 

As elaborated further in Section~\ref{sec:voltage}, the use of the nonlinear relation~\eqref{loadmodel} in $(P1)$ introduces an additional source of non-convexity~\cite{Khodr09},  and would render the DSR formulated in the next section nonconvex. This, in turn, would exacerbate the problem complexity, and would make optimality claims on the obtained topology difficult to establish. Instead, the aim here is to develop a DSR scheme that is computationally efficient yet able to reliably discard lines involving high active power losses. Developing such a scheme is instrumental for network operators to quickly re-shape the distribution grid in case of e.g. abrupt load or generation variations, and to promptly restore the power delivery network-wide after an outage event~\cite{Pilo99}. 

To this end, the powerful  approximate load model derived in~\cite{Bolognani13} is advocated to relate injected currents to complex powers linearly. Specifically, upon defining the vectors $\biota_n :=  [\Re^\cT\{\bi_{n}\}, \Im^\cT\{\bi_{n}\}]^{\cT}$ and $\bsigma_{G,n} := [\Re^\cT\{\bs_{G,n}\}, \Im^\cT\{\bs_{G,n}\}]^\cT$, it follows from~\cite{Bolognani13} that the current injected at node $n$ can be approximated as  
\begin{align}
\label{approx_current}
\biota_n \approx \underbrace{\frac{\sqrt{3}}{V_N} 
\left[ 
\begin{array}{rr}
\Re\{\bPhi_n\} & \Im\{\bPhi_n\} \\
\Im\{\bPhi_n\} & - \Re\{\bPhi_n\} 
\end{array}
\right] (\bsigma_{G,n} - \bsigma_{L,n})}_{:= \bg_n(\bsigma_{G,n}) } 
\end{align}
with $\bPhi_n := \diag(\{e^{\varphi_N^\phi}\}_{\phi \in \cP_n})$, and $\bsigma_{G,n} = \mathbf{0}$ for $n \in \cN \backslash (\{1\} \cup \cS)$.  Then, to account for PQ-loads in (P1), replace $\bi_{n}$ in~\eqref{kcl} with the right hand side of~\eqref{approx_current}. The approximation error incurred by~\eqref{approx_current}  is infinitesimal for large nominal voltages $V_N$ as analytically shown in~\cite{Bolognani13}. For example, for the IEEE 37-node feeder, the average error is just on the order of $0.1$ Ampere (a relative error of less than 2\% - which yields an error on the powers that is on the order of the load prediction error)~\cite{Bolognani13}. The motivation behind~\eqref{approx_current} is twofold: first, using~\eqref{approx_current} in conjunction with compressive sampling methods, offers the advantage of a convex DSR problem for a system featuring PQ loads. Further, from an DSR standpoint, this approximation does not jeopardize the ability of the methods proposed in the ensuing sections to efficiently capture topologies yielding low power losses. Based on the found topology, voltages and currents are as usual fine-tuned in a subsequent stage by employing  more sophisticated techniques such as OPF.

\subsection{\textcolor{black}{Alternative cost functions}} 
\label{sec:cost_functions}

\textcolor{black}{Similar to e.g.,~\cite{Merlin75,Baran89,Shirmohammadi89,Cho93,Khodr09}, the reconfiguration problem~\eqref{recon_nonconvex_current} considers minimizing the overall active power loss. However, alternative
objectives can be pursued as exemplified next.}

\noindent \textcolor{black}{\emph{Cost of supplied power}. Let $c_0$ denote the cost of active power drawn from the distribution substation, and $c_n^\phi$ the one incurred by the use of a DG unit located at phase $\phi$ of node $n \in \cS$. Supposing that DG units operate at unitary power factor, the net network operational cost can be minimized by replacing~\eqref{P1cost} with: 
\begin{align}
C_{op}(\{\bi_{mn}\}) & := c_0 \Big[ \sum_{(m,n) \in \cE} \Re^\cT\{\bi_{mn}\} \Re\{\bZ_{mn}\} \Re\{\bi_{mn}\} \nonumber \\ 
&  \hspace{-1.3cm} + \Im^\cT\{\bi_{mn}\} \Re\{\bZ_{mn}\} \Im\{\bi_{mn}\}  \Big] + \sum_{n \in \cS} \sum_{\cP_n} c_n^\phi \Re\{S_n^\phi\}  \label{Operationalcost} 
\end{align}
where the first term on the right hand side accounts for the cost of power losses on the network. Clearly,~\eqref{Operationalcost} subsumes~\eqref{P1cost}.}

\noindent \textcolor{black}{\emph{Load balancing}. The load balancing index defined in e.g.,~\cite{Baran89,Ke04} for balanced networks can be extended to the unbalanced setup.  To this end, let $\cE_B \subset \cE$ collect branches whose loading condition is to be controlled, and consider adopting the ratio $| I_{mn}^\phi |^2 / I_{mn}^{\textrm{max}}$ as a loading index for conductor $\phi$ of branch $(m,n) \in \cE_B$~\cite{Ke04}. For instance, $\cE_B$ may include transformers, tie lines, or the first segment of (sub-)laterals.  Then, to facilitate a more equitable treatment of branches $\cE_B$ in terms of loading, the following cost can be minimized (see also~\cite{Baran89} and~\cite{Ke04}): 
\begin{align}
C_{bal}(\{\bi_{mn}\}) & :=  \sum_{(m,n) \in \cE_B} \sum_{\phi \in \cP_{mn}} \frac{| I_{mn}^\phi |^2}{I_{mn}^{\textrm{max}}} .\label{Balancing} 
\end{align}
Clearly, a weighted combination of~\eqref{Operationalcost} and~\eqref{Balancing} can also be considered in order to trade off operational costs for system security.   
}

\section{DSR via group-sparsity} 
\label{sec:dsr_sparsity}

Collect first the real and imaginary parts of $\bi_{mn}$ in the vector $\bxi_{mn} := [\Re^\cT\{\bi_{mn}\}, \Im^\cT\{\bi_{mn}\}]^{\cT} \in \mathbb{R}^{2|\cP_{mn}|}$, and define matrices $\bbmZ_{mn} := \bI_{2} \otimes \Re\{\bZ_{mn}\}$, $\bbmA_{mn}^{(m)} := \bI_{2} \otimes \bA_{mn}^{(m)}$, and $\bbmM_{mn}^{\phi} := \bI_{2} \otimes \be_{mn}^\phi (\be_{mn}^\phi)^\cT$, where $\bI_{2}$ is the $2 \times 2$ identity matrix. Key to obtaining a convex re-formulation of $(P1)$ is to notice that the entries of  $\bxi_{mn}$ are \emph{all} zero if line $(m,n) \in \cE_R$ is not used to deliver power to the loads. In compressive sampling, this translates to having the vector $\bxi_R := [\{\bxi_{mn}^\cT | (m,n) \in \cE_R\}]^\cT$ being \emph{group-sparse}~\cite{YuLi06}; meaning that, either $\bxi_{mn} = \mathbf{0}$, or, the elements of $\bxi_{mn}$ may all be nonzero. One powerful way to capitalize on this attribute of currents flowing on lines equipped with switches, consists in augmenting the cost~\eqref{P1cost} with the following sparsity-promoting regularization term~\cite{YuLi06}
\begin{align}
r(\{\bxi_{mn}\}) :=  \sum_{(m,n) \in \cE_{R}} \lambda \,  \| \bxi_{mn} \|_2  \label{Glasso_currents}  
\end{align}
where $\lambda$ is a real positive constant.  Then, using~\eqref{Glasso_currents}, and discarding binary line selection variables, the DSR problem can be re-formulated as: 
\begin{subequations}
\label{recon_convex_current}
\begin{align}
& \hspace{-1.8cm} (P2)  \,  \min_{\{\bxi_{mn}\}, \{\bsigma_{G,n}\}  } \,\, \frac{1}{2} \sum_{(n,m) \in \cE} \bxi_{mn}^\cT \bbmZ_{mn} \bxi_{mn} + r(\{\bxi_{mn}\}) \label{P2cost}  \\
\textrm{subject~to~~}  & \hspace{.5cm}  \bsigma_{G,n}^{\mathrm{min}} \leq \bsigma_{G,n} \leq  \bsigma_{G,n}^{\mathrm{max}} , \, \forall n \in \cS 
\label{P2dg}   \\
& \hspace{-0.9cm}  \bg_n(\bsigma_{G,n}) + \sum_{j \in \cN_{\rightarrow n}} \hspace{-.2cm} \bbmA_{jn}^{(n)} \, \bxi_{jn}  =  \sum_{k \in \cN_{n \rightarrow}}  \hspace{-.2cm} \bbmA_{nk}^{(n)} \, \bxi_{nk} \label{P2kcl} \\
& \hspace{0.2cm} \bxi_{mn}^{\cT} \bbmM_{mn}^{\phi} \bxi_{mn} \leq I_{mn}^{\textrm{max}} , \, \forall \, \phi , \, (m,n) \in \cE  \label{P2current}  
\end{align}
\end{subequations}
where Kirchhoff's current law~\eqref{P2kcl} is enforced at each node $n \in \cN$, and~\eqref{P2dg} are box constraints for the power supplied by controllable DG units. Matrices $\{\bbmZ_{mn}\}$ are positive definite, while matrices $\{\bbmM_{mn}^{\phi}\}$ are symmetric, positive semidefinite, with rank 2. Thus, $(P2)$ is convex, and it can be solved optimally via general-purpose interior point methods. What is more, $(P2)$ can be reformulated as a SOCP, and thus solved using efficient primal-dual interior point methods tailored to SOCP (see e.g.,~\cite{Lobo98}). \textcolor{black}{To this end, one has to introduce the auxiliary variables $\{t_{mn}\}_{(m,n) \in \cE_R}$, replace $r(\{\bxi_{mn}\})$ in~\eqref{P2cost} with the linear term $\sum_{(m,n) \in \cE_R} \lambda t_{mn}$, and add the second order cone constraints $\|\bxi_{mn}\|_2 \leq t_{mn}, (m,n) \in \cE_R$. The optimal topology of the network is obtained by discarding the distribution lines with an associated zero current. That is, 
$\cE^{opt} := \cE \backslash \{(m,n) \in \cE_R : \bxi_{mn}^{opt} = \mathbf{0} \}$.}

\textcolor{black}{When the objective is to minimize the net operational cost, or, promote load balancing, the first term in the cost~\eqref{P2cost} should be replaced by 
 $C_{op}(\{\bi_{mn}\}) = c_0 (\sum_{(n,m) \in \cE} \bxi_{mn}^\cT \bbmZ_{mn} \bxi_{mn}) + \sum_{n \in \cS} \sum_{\cP_n} c_n^\phi (\be_n^\phi (\be_n^\phi)^\cT \sigma_{G,n}$, $C_{bal}(\{\bi_{mn}\}) =  \sum_{(m,n) \in \cE_B} \sum_{\phi \in \cP_{mn}} (1/I_{mn}^{\textrm{max}})  \bxi_{mn}^{\cT} \bbmM_{mn}^{\phi} \bxi_{mn}$, or a combination of the two [cf.~\eqref{Operationalcost},~\eqref{Balancing}]. } 

The role of $\lambda$ in $r(\{\bxi_{mn}\})$ is to control the number of vectors $\{\bxi_{mn}\}_{(m,n) \in \cE_R}$ (and, hence, of currents $\{\bi_{mn}\}_{(m,n) \in \cE_R}$) that are set to zero. When $\lambda = 0$, all branches $\cE_R$ are traversed by a non-zero current. Then, with $\lambda$ increasing, the number of lines where no current is flowing increases~\cite{YuLi06}. This implies that by adjusting $\lambda$ one can obtain either meshed topologies (low values of $\lambda$), weakly-meshed, or even radial systems (high values of $\lambda$). To rigorously show this, results from duality theory~\cite{BertsekasConvexAnOpt} are leveraged next to derive closed form expressions for the optimal line currents.  

Let $\{\bmu_n\}$ and $\{\rho_{mn}^{\phi}\}$ denote the multipliers associated with~\eqref{P2kcl} and~\eqref{P2current}, respectively, and consider the (partial) Lagrangian function of $(P2)$, namely:
\begin{align}
& \cL(\bxi, \bsigma_{G}, \bmu, \brho) := \frac{1}{2} \sum_{(n,m) \in \cE} \bxi_{mn}^\cT \bbmZ_{mn} \bxi_{mn} + r \left(\{\bxi_{mn}\}\right) \nonumber 
\end{align}
\begin{align}
& + \sum_{n = 1}^{N} \bmu_n^{\cT} \left(   \bg_n(\bsigma_{G,n}) + \sum_{j \in \cN_{\rightarrow n}} \hspace{-.2cm} \bbmA_{jn}^{(n)} \, \bxi_{jn}  -  \sum_{k \in \cN_{n \rightarrow}}   \bbmA_{nk}^{(n)} \, \bxi_{nk}  \right) \nonumber \\
& + \sum_{(n,m) \in \cE} \sum_{\phi \in \cP_{mn}} \rho_{mn}^{\phi} \left( \bxi_{mn}^{\cT} \bbmM_{mn}^{\phi} \bxi_{mn} - I_{mn}^{\textrm{max}}   \right)   \label{P2Lagrangian}  
\end{align}
where $\bxi := \{\xi_{mn}\}$, $\bsigma_{G} := \{\bsigma_{G,n}\}$, and likewise $\bmu, \brho$ collect all the dual variables for brevity.
Given~\eqref{P2Lagrangian}, the dual function and the dual problem take the form
\begin{align}
\cD( \bmu, \brho) & := \min_{\bxi, \,\,\bsigma_{G,n}^{\mathrm{min}} \preceq \bsigma_{G,n} \preceq  \bsigma_{G,n}^{\mathrm{max}}} \cL(\bxi, \bsigma_{G}, \bmu, \brho)  \label{dual}  \\
\cD^{opt} & =  \max_{\{\bmu_n \succeq \mathbf{0} \} , \brho}  \cD( \bmu, \brho) \, . \label{opt_duality}  
\end{align}
Since $(P2)$ is convex, if there exists a feasible solution $\bxi, \bsigma_{G}$ such that $\bxi_{mn}^{\cT} \bbmM_{mn}^{\phi} \bxi_{mn}  < I_{mn}^{\textrm{max}}$ for all $(m,n) \in \cE$ and $\bsigma_{G,n}^{\mathrm{min}} \prec  \bsigma_{G,n} \prec  \bsigma_{G,n}^{\mathrm{max}}$ for all $n \in \cS$ (that is, Slater's condition holds), then $(P2)$ has zero duality gap~\cite[Ch.~6]{BertsekasConvexAnOpt}. Suppose that this is the case, and let $\bxi^{opt}, \bsigma_{G}^{opt}$ and $\bmu^{opt}, \brho^{opt}$ denote the optimal primal and dual solutions, respectively. The Lagrangian optimality condition~\cite[Prop.~6.2.5]{BertsekasConvexAnOpt} asserts that $\bxi^{opt}$ and $\bsigma_{G}^{opt}$ are also the minimizers of~\eqref{dual} for $\bmu = \bmu^{opt}$ and $\brho = \brho^{opt}$; that is $\cD^{opt} ( \bmu^{opt}, \brho^{opt}) = \cL(\bxi^{opt}, \bsigma_{G}^{opt}, \bmu^{opt}, \brho^{opt})$. Thus, re-arranging terms of the Lagrangian function in a convenient way, and 
exploiting the decomposability of~\eqref{P2Lagrangian}, it turns out that the optimal currents flowing on the phases of line $(m,n)$ are given as the solution of the sub-problem:
\begin{align}
\hspace{-.2cm} \bxi_{mn}^{opt} = \arg \min_{\bxi_{mn} } \, \frac{1}{2} \bxi_{mn}^\cT \tbmZ_{mn} \bxi_{mn} + \lambda \|\bxi_{mn} \|_2 - \bmu_{mn}^\cT  \bxi_{mn} \hspace{-0.1cm}  \label{Optimal_current}  
\end{align}
where $\lambda = 0$ for lines $(m,n) \in \cE \backslash \cE_R$ (whereas $\lambda > 0$ for all lines $\cE_R$), and 
\begin{align}
\tbmZ_{mn} & := \bbmZ_{mn} + \sum_{\phi \in \cP_{mn}} \rho_{mn}^{\phi, opt} \bbmM_{mn}^{\phi} \\
\bmu_{mn} & :=   \bbmA_{mn}^{(m) \cT} \bmu_m^{opt} -  \bbmA_{mn}^{(n) \cT} \bmu_n^{opt} .
\end{align}
Since $\bbmZ_{mn}$ is positive definite and $\sum_{\phi} \rho_{mn}^{\phi, opt} \bbmM_{mn}^{\phi}$ positive semidefinite, it follows that $\bbmZ_{mn}$ is positive definite and invertible. Thus, based on~\eqref{Optimal_current}, the optimal line currents are obtained next. 

\vspace{.2cm}

\begin{proposition}
\label{prop:current_no_switch}
Per line $(m,n) \in \cE \backslash \cE_R$, the optimal currents $\bxi_{mn}^{opt}$ are given by 
\begin{align}
\bxi_{mn}^{opt} 
= \tbmZ_{mn}^{-1} \, \bmu_{mn}   \, .
\label{Optimal_current_no_switch}  
\end{align}
\end{proposition}
    
\vspace{.2cm}    
    
\begin{proposition}
\label{prop:current_switch}
If $(m,n) \in \cE_R$ is a single-phase distribution line, then the optimal current $\bxi_{mn}^{\phi, opt} = [\Re\{I_{mn}^{\phi}\}, \Im\{I_{mn}^{\phi}\}]^\cT$ on phase $\phi$ is given by the following soft-thresholding vector operation
\begin{align}
\hspace{-.0cm} \bxi_{mn}^{\phi, opt}  =  \frac{\left[\|\bmu_{mn}\|_2 - \lambda \right]_{+}}{\left(\Re\{Z_{mn}^\phi\}  + \rho_{mn}^{\phi, opt} \right) \|\bmu_{mn}\|_2} \bmu_{mn}  
\label{Optimal_current_switch}   
\end{align}
where $[a]_+ := \max\{0,a\}$. For lines $(m,n) \in \cE_R$ that are two- or three-phase, the optimal vector of line currents $\bxi_{mn}^{opt}$ is obtained via the following shrinkage and thresholding vector operation   
\begin{align}
 \bxi_{mn}^{opt}  & =  \eta^{opt} \mathbb{I}_{\left\{\|\bmu_{mn}\|_2 > \lambda \right\}}  \left(\eta^{opt}  \tbmZ_{mn}  + \frac{\lambda^2}{2} \bI_{2|\cP_{mn}|} \right)^{-1} \bmu_{mn}
\label{Optimal_current_switch2}   
\end{align}
where $\eta^{opt} \in \mathbb{R}^+$ is the solution of the scalar optimization problem
\begin{align}
 \min_{\eta \geq 0} \, \eta - \frac{\eta}{2} \bmu_{mn}^\cT \left( \eta \tbmZ_{mn}  + \frac{\lambda^2}{2} \bI_{2|\cP_{mn}|} \right)^{-1} \bmu_{mn}.
 \label{Optimal_current_eta}  
 \end{align}
\end{proposition}

\emph{Proof.} See the Appendix. 

\vspace{.2cm}
  
Some comments are now due in order to interpret the role of the multipliers $\{\bmu_n^{opt}\}$ and $\{\rho_{mn}^{\phi, opt}\}$ in view of Ohm's Law, and to better appreciate the merits of the sparsity-promoting regularization term $r(\{\bxi_{mn}\})$. Notice first that from the complementary slackness condition~\cite[Prop.~6.2.5]{BertsekasConvexAnOpt}, one has that $\rho_{mn}^{\phi, opt} = 0$ whenever the corresponding constraint~\eqref{P2current} is not active. Suppose temporarily that $\rho_{mn}^{\phi, opt} = 0$ for all lines, in which case $\tbmZ_{mn}$ boils down to $\tbmZ_{mn} = \bI_{2} \otimes \Re\{\bZ_{mn}\}$ [cf.~\eqref{P2Lagrangian}].  Since currents and voltages abide by Ohm's Law,~\eqref{Optimal_current_no_switch}--\eqref{Optimal_current_switch} imply that the legitimate unit for $\{\bmu_{mn}\}$ is the volt. In particular, comparing~\eqref{Optimal_current_no_switch} with~\eqref{line_currents} reveals that $\bmu_{mn}^{opt}$ corresponds to the electrical potential difference between two nodes $m$ and $n$ connected by a line with a resistive matrix $\Re\{\bZ_{mn}\}$; that is,  $\bmu_{mn}^{opt} = [\Re^{\cT}\{\Re\{\bZ_{mn}\} \bi_{mn}\}, \Im^{\cT}\{\Re\{\bZ_{mn}\} \bi_{mn}\}]^\cT$. In other words, $\bmu_{mn}$ represents the contribution to the potential difference $\bv_m - \bv_n$ that is due to the resistive part of $\bZ_{mn}$. 

With this  connotation of $\bmu_{mn}$, it follows from~\eqref{Optimal_current_switch} that single-phase lines equipped with a switch can be characterized by a resistance given by $\Re\{Z_{mn}^{\phi}\} \|\bmu_{mn}\|_2/[\|\bmu_{mn}\|_2 - \lambda]_+$. This resistive boost discourages  high currents on line $(m,n)$, something that in compressive sampling is usually referred to as ``shrinkage operation''~\cite{YuLi06,Wiesel11}. Eventually, when $\|\bmu_{mn}\|_2 < \lambda$, the value of this resistance goes to infinity, thus resembling an open switch. 
Notice further that the thresholding operator $\left[\|\bmu_{mn}\|_2 - \lambda \right]_{+}$ naturally suggests the order of magnitude of the parameter $\lambda$ that has to be used to (de)select a line. Although less intuitive, this shrinkage and thresholding operation effected through  $\lambda \|\bxi_{mn}\|_2$ can be noticed also in~\eqref{Optimal_current_switch2} for lines with two and three phases.  Here, the design variable $\eta$ is expressed in watt.   
Finally, variable $\rho_{mn}^{\phi, opt}$ can be interpreted as an additional resistance added to the conductor $\phi$ of line $(m,n)$ when the currents reaches its maximum allowable value $I_{mn}^{\mathrm{max}}$. In principle, since this extra resistive value introduce an additional power loss, higher values of $I_{mn}^{\phi}$ are discouraged.    

 
Finally, notice that in order to encourage the use of specific lines $(m,n) \in \cE_R$, the regularization function~\eqref{Glasso_currents} can be replaced by its weighted counterpart $r^{\prime} (\{\bxi_{mn}\}) :=  \sum_{(m,n) \in \cE_R} \lambda_{mn} \| \bxi_{mn} \|_2$, with $\lambda_{mn} \geq 0$ for all $(m,n)$. For example, if the use of a line $(j,k)$ is inadvisable due to e.g., ongoing maintenance or security concerns, a higher associated weight $\lambda_{jk} > \lambda_{mn}$ should be selected.  

\vspace{.1cm}

\textcolor{black}{
\emph{Remark 1~(distribution transformers).} Efficiency of distribution transformers is defined as\footnote{\textcolor{black}{See e.g., U.S. Code of Federal Regulations, Title 10 (Department of Energy), Chapter II, Part 431, Subpart K, Edition Jan. 2013. [Online] Available at: \texttt{http://www.gpo.gov}; and, T. R. Blackburn, ``Distribution transformers: Proposal to increase MEPS levels,'' Technical report, Oct. 2007.  [Online] Available at: \texttt{http://www.energyrating.gov.au.}}}
\begin{align}
\eta_n^\phi := \frac{\bar{P}_{L,n}^\phi}{\bar{P}_{L,n}^\phi + L_{c,n}^\phi + L_{nc,n}^\phi}
 \label{efficiencyTransformer}  
 \end{align}
 where $\bar{P}_{L,n}^\phi$ is the power demanded at the secondary of the transformer; $L_{c,n}^\phi > 0$ is the no-load core loss; and, $L_{nc,n}^\phi > 0$ represents the load loss, which emerges from ohmic losses in the windings, stray losses, core clamps, magnetic shields, and other conducting parts. Values for $\eta_n^\phi$ are higher 
than $98\%$ and, therefore, the voltamperes of the load (secondary winding) and the source (primary winding) are typically assumed to coincide (see also~\cite[Ch.~8]{Kerstingbook}). This explains why in prior works on system reconfiguration~\cite{Merlin75,Baran89,Shirmohammadi89,Cho93,Schmidt95,Gomes06,Huang02,Khodr09,Moradzadeh12} and OPF~\cite{Tse12,Dallanese-TSG13}, the distribution transformer losses were not accounted for explicitly.  }
 
\textcolor{black}{The no-load core loss $L_{c,n}^\phi$ is a transformer-specific constant, evaluated at a reference temperature. On the other hand, the load loss varies with the primary and secondary currents, and can be approximated as $L_{c,n}^\phi = | I_{n}^{\phi} |^2 R_{n}^\phi$, where $I_n^\phi$ is taken to be the current on the primary and $R_{n}^\phi > 0$ is a resistive coefficient. Specifically, $R_{n}^\phi > 0$ is computed based on the transformer turn ratio, the dc primary and secondary resistances, and given temperature correction factors. From~\eqref{efficiencyTransformer}, it follows that the active power required at the primary of the transformer can be expressed as 
\begin{align}
P_{L,n}^\phi = \frac{\bar{P}_{L,n}^\phi}{\eta_n^\phi} = \bar{P}_{L,n}^\phi + L_{c,n}^\phi + |I_{n}^\phi|^2 R_{n}^\phi \, .
\label{powePrimaryTransformer}  
 \end{align} 
Thus, to account for  distribution transformer losses, one can readily replace $\Re\{S_{L,n}^\phi\}$ with $\bar{P}_{L,n}^\phi + L_{c,n}^\phi + |I_{n}^\phi|^2 R_{n}^\phi$ in Section~\ref{sec:load_approximation}.
} 

\vspace{.1cm}

\textcolor{black}{
\emph{Remark 2~(computational complexity).} Solving the reconfiguration problem using e.g., branch and bound techniques~\cite{Cho93} or other heuristic schemes~\cite{Baran89,Shirmohammadi89,Schmidt95,Gomes06} incurs higher computational burden than $(P2)$. The reason is twofold: \emph{i)} heuristics in~\cite{Baran89,Shirmohammadi89,Cho93,Schmidt95,Gomes06} are iterative methods that require testing \emph{multiple} combinations of the binary variables, and a nonlinear power flow problem must be solved for each combination; \emph{ii)} power flow problems can be solved by using either off-the-shelf solvers for nonlinear programs~\cite{Khodr09}, or, semidefinite relaxation schemes as in~\cite{Tse12,Dallanese-TSG13}. Either way, these methods incur complexity higher than that of a single SOCP. In fact, letting $r$ denote the total number of optimization variables, and $\epsilon > 0$ a prescribed solution accuracy, the worst-case complexity on the order of $\mathcal{O}(r^3 \log(1/\epsilon))$ for SOCPs, $\mathcal{O}(r^{4.5})\log(1/\epsilon))$ for SDPs~\cite{Nesterov94}, and typically even higher for solvers for nonlinear programs~\cite{PaudyalyISGT}.  Notice also that the complexity of SOCPs scales better as the system size $r$ increases. }

\section{Accounting for voltage constraints}
\label{sec:voltage}

Similar to various DSR renditions, the objective of $(P2)$ is to obtain a topology that is likely to yield  the lowest power losses for a given predicted load profile. Based on the resultant optimal configuration $(\cN,\cE^{opt})$, voltages and currents are fine-tuned in a subsequent optimization stage where more sophisticated techniques such as OPF are employed (see e.g.,~\cite{Dallanese-TSG13} and references therein). In some cases however, it may be desirable to introduce voltage regulation-related constraints in order to avoid network configurations that can potentially yield infeasible OPF solutions (meaning a set of voltages returned by the OPF solver not within prescribed minimum and maximum utilization limits). 

To effect voltage regulation, consider introducing a constraint $V_n^{\phi} \in \cB_n^{\phi}$ per node $n \in \cN \backslash \{1\}$ and phase $\phi \in \cP_n$, where $\cB_n^{\phi}$ is a given closed set collecting the admissible voltages (see e.g.,~\cite{Khodr09, Kerstingbook}); set $\cB_n^{\phi}$ will be exemplified in the ensuing Sections~\ref{sec:voltagebox} and~\ref{sec:voltagemagnitude}. One way to enforce these constraints is to let the voltages become optimization variables, and formulate a joint DSR and OPF problem as in, e.g.~\cite{Khodr09}. However, it is not convenient here to have voltages as optimization variables because:

 \emph{i)} equality~\eqref{loadmodel} is nonconvex (as in OPF problems); and, 
 
 \emph{ii)} when binary variables are used to model the states of switches~\cite{Cho93,Schmidt95,Gomes06,Huang02,Khodr09,Moradzadeh12}, and constraints~\eqref{P1close} are employed~\cite{Cho93,Khodr09,Moradzadeh12}, the solver would set $\bi_{mn}^{opt} = \mathbf{0}$ for lines with $x_{mn} = 0$. However, it is clear from~\eqref{line_currents} that imposing $\bi_{mn}^{opt} = \mathbf{0}$ requires equating voltages at the two end points of line $(m,n)$; that is, $[\bv_m]_{\cP_{mn}} = [\bv_n]_{\cP_{mn}}$. This artifact renders the joint DSR and OPF problem infeasible in various practical cases. Consider for example a network with $5$ nodes, lines $\cE = \{(1,2), (2,3),(2,4),(2,5),(3,5),(4,5)\}$, and switches in $\cE_{R} = \{(2,3), (2,4), (2,5)\}$. Suppose that only one switch must be closed in order to obtain a radial network, and this switch is the one on line $(2,3)$. However, by setting $\bv_2 = \bv_4 = \bv_5$, the load demands at nodes 4 and 5 would not be satisfied, since no power is delivered at these nodes. 
 
 One approach to resolving this issue is to discard~\eqref{P1close} and replace $\{\bZ_{mn}\}_{(m,n) \in \cE_R}$ with $\{\bZ_{mn} x_{mn} \}_{(m,n) \in \cE_R}$. However, the resultant formulation yields  a challenging bilinear problem with integer variables. Aiming at a computationally affordable DSR scheme, voltages are treated here as \emph{latent} problem variables as shown next.
 
Recall from Section~\ref{sec:load_approximation} that voltages at the substation  $\bnu_1 :=  [\Re^\cT\{\bv_1\}, \Im^\cT\{\bv_1\}]^{\cT}$ are typically taken as a reference for the phasorial representation~\cite{Kerstingbook}. Consider the network in Fig.~\ref{fig:F_feeder37}. Given $\bnu_1$, and assuming that the arc $(1,2)$ originates at node $1$ and ends at $2$, voltages at node $2$ can be expressed as [cf.~\eqref{line_currents}]
\begin{align}
\label{voltage_drop}
\bnu_2 = \bnu_1 -  \underbrace{
\left[ 
\begin{array}{rr}
\Re\{\bZ_{12}\} & -\Im\{\bZ_{12}\} \\
\Im\{\bZ_{12}\} & \Re\{\bZ_{12}\} 
\end{array}
\right]}_{:= \bPsi_{12} } \bxi_{12} \, .
\end{align}
Likewise, if the sequence of nodes $\cW_{1\rightarrow n} := \{1, \ldots, n\}$ forms an undirected path $1 \rightarrow 2 \rightarrow ... \rightarrow n$ from the substation to node $n$, and \emph{none} of the lines $(m,n): m,n \in \cW _{1\rightarrow n}$ is equipped with switches, then $\bnu_n$ can be expressed as $\bnu_n = \bnu_1 + \sum_{(m,n): m,n \in \cW _{1\rightarrow n} }  \alpha_{mn}  \bPsi_{mn}  \bxi_{mn}$, where $\alpha_{mn} = 1$ if the path traverses the directed edge $(m,n)$ (which goes from $m$ to $n$) in the opposite direction, and $\alpha_{mn} = -1$ otherwise. 
A similar approach was taken in~\cite{Cho93}.
Based on these relations, voltage regulation can be readily effected by adding to $(P2)$ the following constraint per node $n$:
\begin{align}
\label{voltage_constr}
\bnu_1 + \hspace{-.3cm}  \sum_{(m,n): m,n \in \cW _{1\rightarrow n} }  \hspace{-.5cm} \alpha_{mn}  \bPsi_{mn}  \bxi_{mn} \in \cB_n \, .
\end{align}
When switches are densely deployed, there may not exist an undirected path connecting the substation to a node $n$ that includes
 only lines in $\cE \backslash \cE_R$.  If this is the case, the substation must be replaced by another point of reference; that is, a node where the value (or an approximate value) of the voltages can be unequivocally determined. From the model set forth in Section~\ref{sec:load_approximation} (see also~\cite{Bolognani13}), it can be noticed that approximate values of the voltages are readily available for nodes with a nonzero load demand.  Hereafter, the reference node will be generically denoted by $n_{ref}$.  Two possible choices for the set $\cB_n$ are presented in the ensuing subsections.

\subsection{Box constraints}
\label{sec:voltagebox}

Let $\hat{V}_{n,\Re}^{\phi}$ and $\check{V}_{n,\Re}^{\phi}$ be upper and lower bounds, respectively, on the real part of voltage $V_n^{\phi}$. Likewise, let $\hat{V}_{n,\Im}^{\phi}$ and $\check{V}_{n,\Im}^{\phi}$ denote the counterparts for $\Im\{V_n^{\phi}\}$. Then, upon collecting these quantities for all phases $\phi \in \cP_n$ in the vectors $\hat{\bnu}_{n} := [\{\hat{V}_{n,\Re}^{\phi}\}, \{\hat{V}_{n,\Im}^{\phi}\}]^\cT$ and $\check{\bnu}_{n} := [\{\check{V}_{n,\Re}^{\phi}\}, \{\check{V}_{n,\Im}^{\phi}\}]^\cT$, $(P2)$ can be readily reformulated as follows:
\begin{align}
& \hspace{-1.8cm} (P2^\prime)  \,  \min_{\{\bxi_{mn}\}, \{\bsigma_{G,n}\}  } \,\, \frac{1}{2} \sum_{(n,m) \in \cE} \bxi_{mn}^\cT \bbmZ_{mn} \bxi_{mn} + r(\{\bxi_{mn}\}) \nonumber  \\
\textrm{subject~to~~}  & \eqref{P2dg}, \eqref{P2kcl}, \eqref{P2current}, ~\mathrm{and} \nonumber \\
& \hspace{-1.2cm}  \check{\bnu}_{n}  \preceq \bnu_{n_{ref}} + \hspace{-.3cm}  \sum_{(m,n): m,n \in \cW _{n_{ref} \rightarrow n} }  \hspace{-.5cm} \alpha_{mn}  \bPsi_{mn}  \bxi_{mn} \preceq \hat{\bnu}_{n} \, . \label{boxconstraintV}
\end{align}
Problem $(P2^\prime)$ is convex, and thus efficiently (and optimally) solved via either general-purpose interior point methods, or, primal-dual schemes tailored to SOCP~\cite{Lobo98}. Through~\eqref{boxconstraintV} it is possible to constrain both the magnitude of $V_n^\phi$ and its deviation from the nominal phase $\angle V_{n_{ref}}^\phi$.  Further, introducing these additional constraints does not alter the expressions for the optimal line currents provided in Propositions~\ref{prop:current_no_switch} and~\ref{prop:current_switch}. 

\subsection{Nonconvex constraints on voltage magnitudes}
\label{sec:voltagemagnitude}
 
In traditional OPF approaches, it is usual to consider lower and upper bounds on the voltage magnitudes, which habitually coincide with minimum and maximum utilization and service voltage levels~\cite{Kerstingbook}. Consider then introducing in $(P2)$ the constraints $\hat{V}_n \leq |V_n^{\phi}| \leq \check{V}_n$, per node $n$ and phase $\phi \in \cP_n$. This leads to the following DSR reformulation: 
 \begin{align}
& \hspace{-1.8cm} (P2^{\prime \prime})  \,  \min_{\{\xi_{mn}\}, \{\bsigma_{G,n}\}  } \,\, \frac{1}{2} \sum_{(n,m) \in \cE} \bxi_{mn}^\cT \bbmZ_{mn} \bxi_{mn} + r(\{\bxi_{mn}\}) \nonumber  \\
\textrm{subject~to~~}  & \eqref{P2dg}, \eqref{P2kcl}, \eqref{P2current}, ~\mathrm{and} \nonumber \\
& \hspace{-1.5cm}  \hat{V}_n^2 \leq  \left(   \bnu_{n_{ref}} +  \sum_{m,n} \alpha_{mn}  \bPsi_{mn}  \bxi_{mn}  \right)^\cT \bbmQ_n^{\phi}  \nonumber 
\end{align}
 \begin{align}
& \hspace{-0.8cm} \times \left(   \bnu_{n_{ref}} +  \sum_{m,n} \alpha_{mn}  \bPsi_{mn}  \bxi_{mn}  \right)  \leq  \check{V}_n^2, \phi \in \cP_n \label{constraintmagV}
\end{align}
 where $\bbmQ_n^{\phi} := \bI_2 \otimes \diag(\be_n^{\phi})$. Unfortunately, the constraint $ - \left(   \bnu_{n_{ref}} + \bPsi  \bxi  \right)^\cT \bbmQ_n^{\phi}\left(   \bnu_{n_{ref}} +  \bPsi \bxi \right) + \hat{V}_n^2  \leq 0$, where $\bPsi$ and  $\bxi$ are defined in the obvious way, is nonconvex since the function on the left hand side is concave. However, this source of non-convexity can be efficiently addressed by resorting to convex approximation techniques. Among candidate methods, the successive convex approximation (SCA) approach proposed in~\cite{MaWr78} is well suited for the problem at hand
because it guarantees first-order Karush--Kuhn--Tucker (KKT) \emph{optimality} under mild regularity conditions. 

To briefly illustrate the general SCA method, consider an optimization problem
\begin{subequations}
\label{nonconvex_example}
\begin{align}
& \min_{\bxi \in \cI} f_0(\bxi) \\
\textrm{subject to } & f_k(\bxi) \leq 0,\quad k=1,2,\ldots,K
\end{align}
\end{subequations}
where $f_0(\bxi)$ is convex and differentiable, $f_k(\bxi)$,
$k=1,\ldots,K$, are differentiable functions, and the feasible
region $\cF := \{\bxi \in \cI|f_k(\bxi) \leq 0, k=1,\ldots,K\}$ is
compact. Then, starting from a feasible point $\bxi^{(0)} \in \cF$, a
series of approximate problems can be solved to locate a 
KKT optimal point of
the original (non-convex) problem. For each $k = 1,\ldots,K$, let
$\tilde f_k(\bxi;\bxi^{(j)})$ denote the surrogate function for
$f_k(\bxi)$, which may depend on the solution $\bxi^{(j)}$ to the
problem in the $(j-1)$-st iteration. The approximate
problem to solve in iteration $j$ is
\begin{subequations}
\begin{align}
& \min_{\bxi \in \cI} f_0(\bxi) \\
\textrm{subject to } & \tilde f_k(\bxi; \bxi^{(j)}) \leq 0, \quad k=1,2,\ldots,K
\end{align}
\end{subequations}
whose feasible region is denoted as $\cF^{(j)}$. Provided that
$\tilde f_k(\bxi;\bxi^{(j)})$ satisfies the following conditions
c1)--c3) for each $k = 1,\ldots,K$, the series of solutions $\bxi^{(j)}$, $j=1,2,\ldots$,
to the approximate problems converges to the KKT point of the
original problem~\eqref{nonconvex_example}:
\begin{enumerate}
\item[c1)] $f_k(\bxi) \leq \tilde f_k(\bxi;\bxi^{(j)}),\quad \forall \bxi \in \cF^{(j)}$

\item[c2)] $f_k(\bxi^{(j)}) = \tilde f_k(\bxi^{(j)}; \bxi^{(j)})$

\item[c3)] $\nabla f_k(\bxi^{(j)}) = \nabla \tilde f_k(\bxi^{(j)};\bxi^{(j)})$.
\end{enumerate}

In order to apply the SCA method to $(P2^{\prime \prime})$, an appropriate surrogate constraint for the nonconvex lower bound in~\eqref{constraintmagV} needs to be determined. Since, $ - \left(   \bnu_{n_{ref}} + \bPsi  \bxi  \right)^\cT \bbmQ_n^{\phi}\left(   \bnu_{n_{ref}} +  \bPsi \bxi \right)$ is quadratic and concave, a linear approximation around the feasible point $\{\bxi_{mn}^{(j)}\}$ satisfying c1)--c3) can be readily found. In fact, after standard manipulations, it follows that the convex problem to be solved at iteration $j$ of the SCA algorithm is:
\begin{subequations}
 \begin{align}
& \hspace{-1.8cm} (P2^{\prime \prime (j)})  \,  \min_{\{\xi_{mn}\}, \{\bsigma_{G,n}\}  } \,\, \frac{1}{2} \sum_{(n,m) \in \cE} \bxi_{mn}^\cT \bbmZ_{mn} \bxi_{mn} + r(\{\bxi_{mn}\}) \nonumber  \\
\textrm{subject~to~~}  & \eqref{P2dg}, \eqref{P2kcl}, \eqref{P2current}, ~\mathrm{and} \nonumber \\
& \hspace{-0.9cm} (\ba_n^{\phi}(\bxi^{(j)}))^{\cT} (\bxi^{(j)} - \bxi) - b_n^{\phi}(\bxi^{(j)}) + \hat{V}_n^2 \leq 0  \label{constraintmagV_lower}  \\
& \hspace{-0.3cm}   \left(   \bnu_{n_{ref}} + \bPsi \bxi  \right)^\cT \bbmQ_n^{\phi} \left(   \bnu_{n_{ref}} +  \bPsi \bxi  \right)  \leq  \check{V}_n^2  \label{constraintmagV_upper} 
\end{align}
with 
\begin{align}
\hspace{-0.2cm} \ba_n^{\phi}(\bxi^{(j)}) & := 2 \bPsi^{\cT} \bbmQ_n^{\phi} \left(  \bnu_{n_{ref}} + \bPsi \bxi^{(j)}  \right) \\
\hspace{-0.2cm} b_n^{\phi}(\bxi^{(j)}) & :=  \left(   \bnu_{n_{ref}} + \bPsi \bxi^{(j)}  \right)^\cT \bbmQ_n^{\phi} \left(   \bnu_{n_{ref}} +  \bPsi \bxi^{(j)}  \right)  .
\end{align}
\end{subequations}
The resultant SCA algorithm amounts to solving a sequence of convex problems $(P2^{\prime \prime (0)}) \rightarrow ... \rightarrow  (P2^{\prime \prime (j)}) \rightarrow ... $, where the feasible point $\bxi^{(j)}$ is taken to be the solution of $(P2^{\prime \prime (j-1)})$ (and $\bxi^{(0)}$ is an arbitrary initial feasible point). Convergence of this scheme is formalized next (see also~\cite{MaWr78}). 

\vspace{0.1cm}

\begin{proposition}
\label{prop:kktP2pp}
Given an initial feasible point $\bxi^{(0)}$, the iterates $\bxi^{(j)}$ generated by $(P2^{\prime \prime (j)})$, $j = 0,1,2, \ldots$, converge to a KKT solution of the nonconvex DSR problem $(P2^{\prime \prime})$.
\end{proposition}

\vspace{0.1cm}

In the ensuing section, the proposed DSR algorithms are tested \textcolor{black}{on three different distribution systems}.
 
\section{Numerical Experiments}
\label{sec:Simulation}

The optimization package \texttt{CVX}\footnote{[Online] Available: \texttt{http://cvxr.com/cvx/}}, along with the interior-point based solver \texttt{SeDuMi}\footnote{[Online] Available: \texttt{http://sedumi.ie.lehigh.edu/}} are employed to solve the DSR problems in \texttt{MATLAB}.  A machine with Intel Core i7-2600 CPU @ 3.40GHz is used.

\textcolor{black}{To implement the proposed reconfiguration strategy for distribution 
systems, the network operator requires: \emph{i)} the phase impedance matrix of each branch; \emph{ii)} the set $\cE_R$ of switches that can be controlled; \emph{iii)} the sets $\{\cN_{n \rightarrow}, \cN_{\rightarrow n}\}$ to be used in~\eqref{kcl}; \emph{iv)} the instantaneous demanded powers; \emph{v)} the minimum and maximum power that DG units can supply; and, \emph{vi)} the parameter $\lambda$ (or, parameters $\lambda_{mn}$ for a weighted regularization function).  
With these data, $(P2)$ can be formulated by following the definitions summarized in Table~\ref{tab:nomenclature}, and solved using standard solvers such as~\texttt{SeDuMi} and~\texttt{CPLEX}\footnote{See \texttt{http://www-01.ibm.com/software/commerce}}, or customized primal-dual algorithms. There are (at least) two viable ways to select $\lambda$. First, since loads typically manifest daily/seasonal patterns, $\lambda$ can be selected based on historical reconfiguration results. Alternatively, when a primal-dual algorithm is employed to solve $(P2)$, $\lambda$ can be adjusted during the iterations of this scheme to obtain a topology with a given number of closes switches. }

\subsection{IEEE 37-node feeder}
\label{sec:feeder1}

\begin{figure}[t]
\begin{center}
\includegraphics[width=0.30\textwidth]{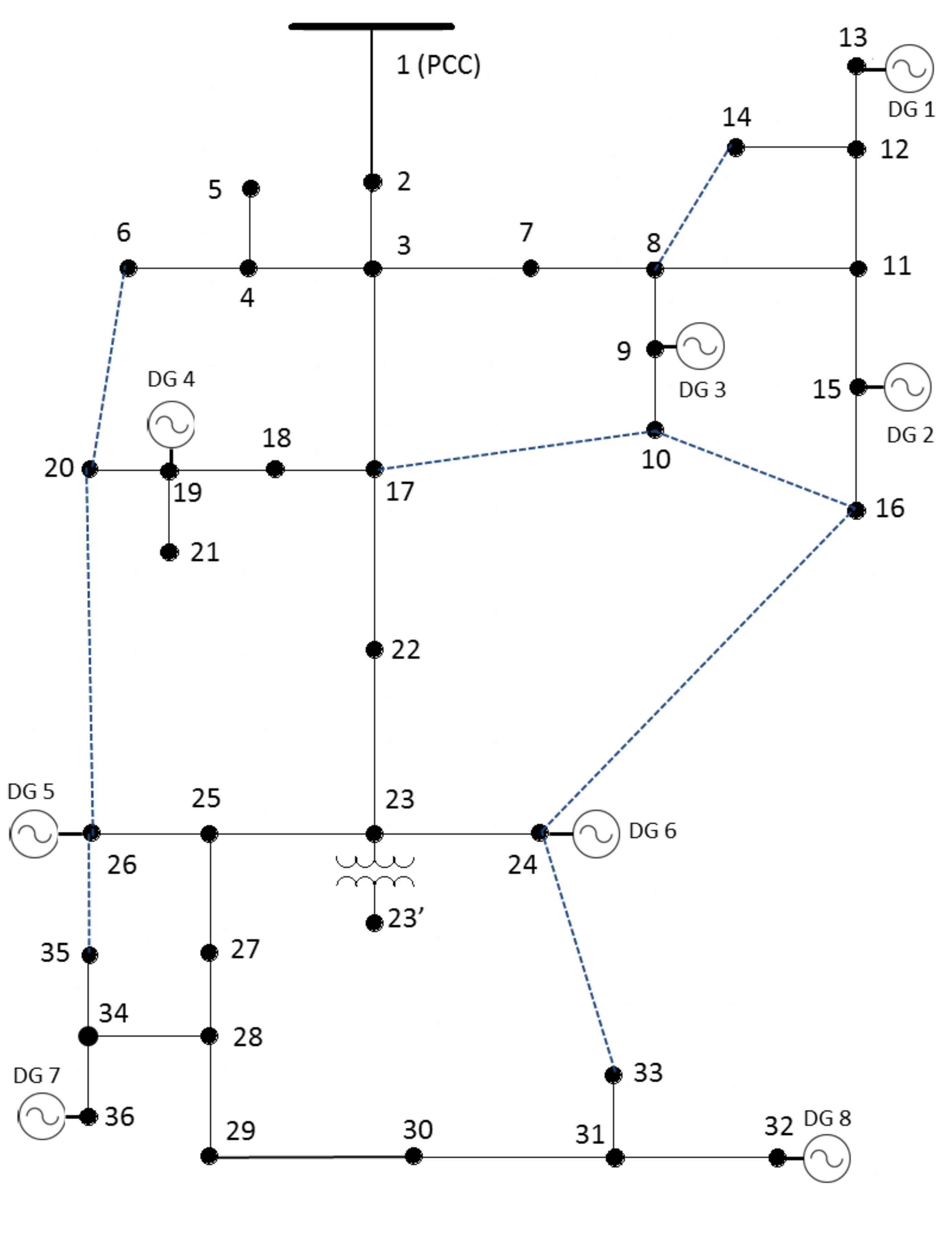}
\vspace{-.5cm}
\caption{Modified IEEE 37-bus test feeder. }
\label{fig:F_feeder37}
\vspace{-.5cm}
\end{center}
\end{figure}

Consider the IEEE 37-node feeder~\cite{testfeeder}, which is a portion of the three-phase $4.8$ kV distribution network located in California. Compared to the scheme in~\cite{testfeeder}, 8 additional three-phase lines equipped with sectionalizing switches are considered, as shown in Fig~\ref{fig:F_feeder37}. The parameters of these additional lines are listed in Table~\ref{tab:lines}, where the admittance matrices corresponding to the configuration indexes 723 and 724 can be found in~\cite{testfeeder}. 
The line impedance matrices for the original lines are computed as specified in~\cite{testfeeder}.
\textcolor{black}{All the demanded complex powers are also listed in~\cite{testfeeder}, and they are modeled as spot constant-PQ loads; thus,~\eqref{approx_current} is used to approximate the injected currents. Distribution transformer losses are neglected}. A balanced load of 85 kW and 40 kVAr per phase is added to node 23 to represent an additional residential demand. Further, controllable DG units are located at nodes $9, 13, 15, 19, 26, 32,$ and $36$; they operate at unity power factor; and, they can supply a maximum power of $50$ kW per phase. Two different setups are considered, depending on the number and locations of switches:
\begin{enumerate}
	\item Test 1: $\cE_{R,1} = \{(8,14), (6,20), (10,16), (20,26),$ $(16,24), (10,17), (24,33), (26,35)\}$; and, 
	\item Test 2: $\cE_{R,2} = \cE_{R,1} \cup \{(17, 22), (23,24), (23,25), $ $(29,30)\}$.
\end{enumerate}
The voltage at the distribution substation is set to $\bv_1 = [1 \angle 0^\circ, 1 \angle -120^\circ, 1 \angle 120^\circ]^\cT$ pu. Finally, box constraints on the voltages are considered; specifically, the lower and upper bounds of the real and imaginary parts of the voltages are such that $\Re\{V_m^{\phi}\} \in [\Re\{V_1^{\phi}\} - 0.0354, \Re\{V_1^{\phi}\} + 0.0354]$ pu and $\Im\{V_m^{\phi}\} \in [\Im\{V_1^{\phi}\} - 0.0354, \Im\{V_1^{\phi}\} + 0.0354]$ pu. This translates to having the magnitude of the voltages in the range $|V_{n}^{\phi}| \in [0.95,1.05]$ pu. 

\begin{table}[h]
\caption{Additional lines in the modified IEEE 37-node feeder}
\begin{center}
\begin{tabular}{|c|cc||c|cc}
Line & Conf. & Length (ft) & Line & Conf. & Length (ft) \\
\hline
(8,14) & 723 & 1144 & (16,24) & 724 & 1580\\
(6,20) & 724 & 1320  & (10,17) & 724 & 1137  \\
(10,16) & 724 & 847 & (24,33) & 724 & 1315 \\
(20,26) & 724 & 815 & (26,35) & 724 & 377 \\
\end{tabular}
\end{center}
\label{tab:lines}
\end{table}%

\begin{figure}
	\centering
  \subfigure[Test 1.]{\includegraphics[width=0.43\textwidth]{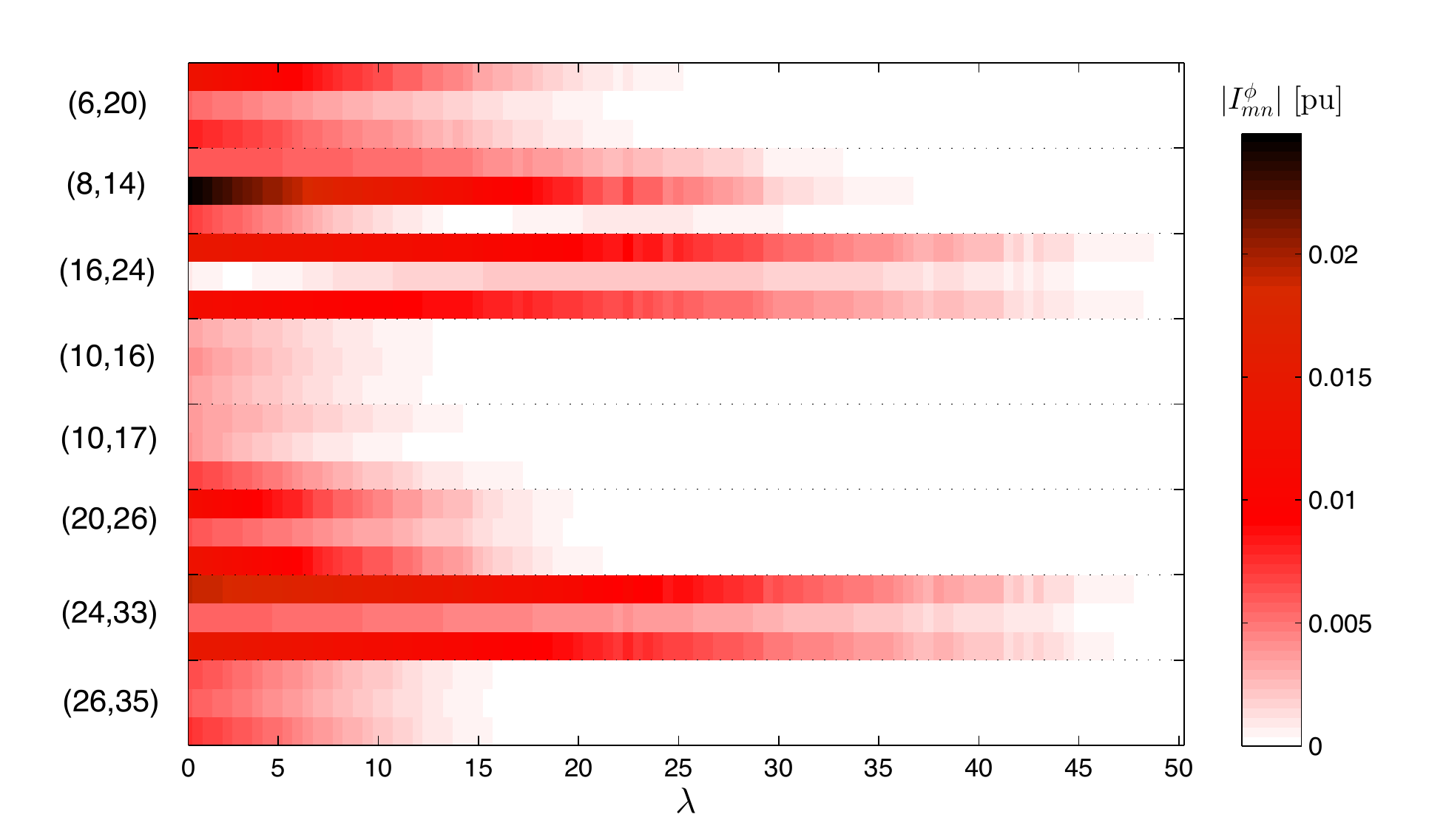}}
  \subfigure[Test 2.]{\includegraphics[width=0.43\textwidth]{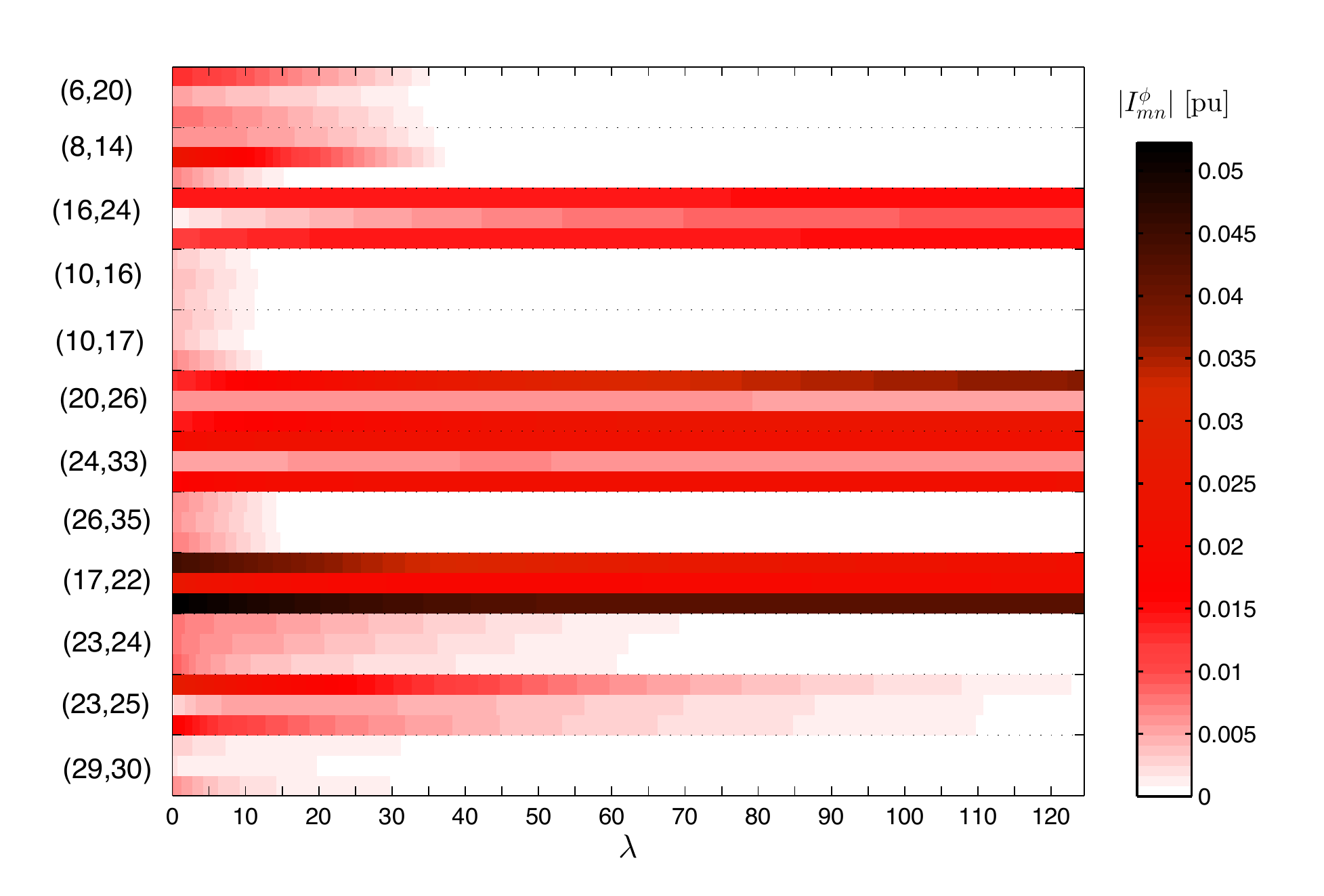}}
  \caption{Absolute value of the currents on lines $\cE_R$, for different $\lambda$.}
  \label{fig:currents}
  \vspace{-.3cm}
\end{figure}

Fig.~\ref{fig:currents} depicts the magnitude of currents flowing on lines equipped with switches, for different values of the sparsity-promoting parameter $\lambda$. Three rows per line are reported, where the first row corresponds to phase ``a,'' and the third one to phase ``c'' (all lines are three-phase). The \textcolor{black}{current magnitude} is color-coded, where white represents a zero current (that is, an open switch), while red hues are used to capture different values for $|I_n^{\phi}| > 0$, in pu.  In the upper plot (a), it can be clearly seen that the number of open switches increases as $\lambda$ increases, and the original tree topology described in~\cite{testfeeder} is obtained at saturation. Interestingly, the first lines that are discarded are $(10, 16)$ and $(10,17)$, which implies that the majority of the power supplied to that part of the network comes from the DG units. As expected, similar trends are evidenced in Fig.~\ref{fig:currents}(b) for Test 2, thus further corroborating the merits of the proposed method based on sparsity-promoting regularizations. It can be noticed that at saturation only four switches are left closed and, interestingly,  the resultant tree topology is different than the original one in~\cite{testfeeder}. In fact, the switches  in $(23,24)$, $(23,25)$ are open, while lines $(16,24)$, $(20,26)$, and $(24,33)$ are used to deliver power to the loads.

\begin{figure}[t]
\begin{center}
\includegraphics[width=0.50\textwidth]{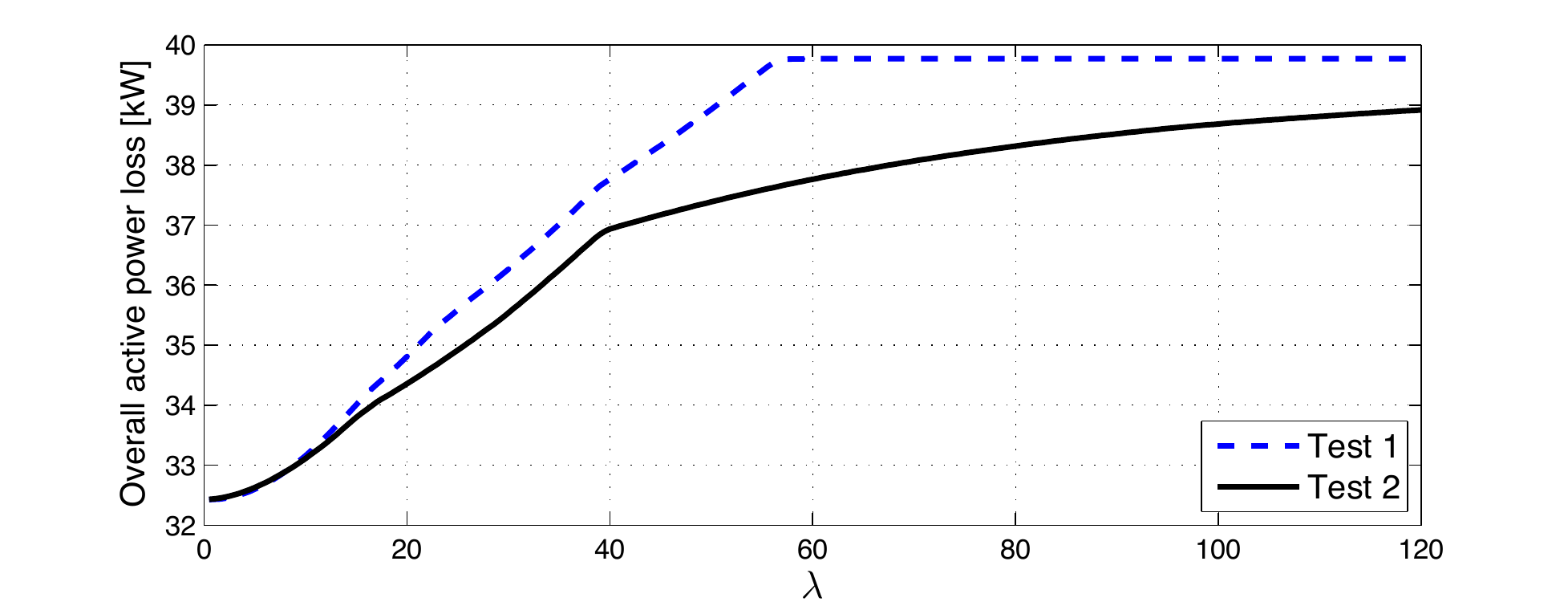}
\vspace{-.3cm}
\caption{Overall active power loss [kW].}
\label{fig:power_loss}
\end{center}
\vspace{-.3cm}
\end{figure}

Fig.~\ref{fig:power_loss} reports the expected active power loss as a function of $\lambda$, for both cases. Notice first that power losses are in general lower in Test 2. Then, it can be clearly seen that the power loss monotonically increases as $\lambda$ increases (that is, with the number of open switches). This motivates augmenting the cost of $(P2)$ with a term that accounts for the possible maintenance costs of lines and switches, in order to find a possible trade-off between active power loss and number of utilized lines. 

Since the linear relation~\eqref{approx_current} introduces an approximation error when PQ loads are present, the average deviation from the nominal loads is quantified next. Specifically, the deviation for the real power is defined as $\Delta P = (1/ \sum_n |\cP_n^L |) \sum_n \sum_\phi |P_{n}^{\phi} - \hat{P}_{n}^{\phi}|$, where $\cP_n^L$ is a set collecting the phases at node $n$ with a non-zero load; $\hat{P}_{L,m}^{\phi}$ is the output of the DSR scheme; and $P_{n}^{\phi} = P_{G,n}^{\phi} - 
P_{L,n}^{\phi}$ is the true real power that would be obtained by considering the exact nonlinear relation~\eqref{loadmodel}. The counterpart $\Delta Q$ is defined in a similar way. These deviations yield only $\Delta P =  1.49$ kW and $\Delta Q = 0.92$ kVAr for Test 1, and $\Delta P =  1.35$ kW and $\Delta Q = 0.85$ kVAr for Test 2. 

\begin{table*}[t]
\caption{\textcolor{black}{Modified IEEE 37-node feeder, test 2, radial system: comparison with other methods}}
\begin{center}
\textcolor{black}{\begin{tabular}{|c||c|c|c|c|}
Method & Open switches &  $P_{\mathrm{loss}}$ [kW] & Time for solver [s] &  Overall time [s] \\
\hline
Proposed  & (6,20), (8,14), (10,16), (10,17), (26,35), (23,24), (23,25), (29,30)  & 41.45 & 0.3 & 2.8 \\
Exhaustive search & (6,20), (8,14), (10,16), (10,17), (26,35), (20,26), (16,24), (29,30) & 41.09 & 235 & 2871 \\
\cite{Shirmohammadi89} (w/ OPF of~\cite{Dallanese-TSG13}) & (6,20), (8,14), (10,16), (10,17), (26,35), (20,26), (16,24), (29,30) & 41.09 & 72.2 & 672 \\
\cite{Gomes06} (w/ OPF of~\cite{Dallanese-TSG13})  & (6,20), (8,14), (10,16), (10,17), (26,35), (20,26), (16,24), (29,30) & 41.09 & 72.2 &  673 \\
Actual network~\cite{testfeeder} & (6,20), (8,14), (10,16), (10,17), (26,35), (20,26), (16,24), (24,31) & 42.37 & --& --
\end{tabular}}
\end{center}
\label{tab:case1_comparison}
\end{table*}%

\begin{table*}[t]
\caption{\textcolor{black}{33-node network in~\cite{Baran89}: comparison with other methods}}
\begin{center}
\textcolor{black}{\begin{tabular}{|c||c|c|c|c|}
Method & Open switches &  $P_{\mathrm{loss}}$ [kW] & Time for solver [s] &  Overall time [s] \\
\hline
Proposed  & (6,7), (9,10), (13,14), (31,32), (24,28)  & 140.28 & 0.3 & 2.7 \\
Proposed w/ weights & (6,7), (8,9), (13,14), (31,32), (24,28)  & 139.56 & 0.3 & 2.7 \\
Exhaust. search~\cite{Morton00} & (6,7), (8,9), (13,14), (31,32), (24,28)  & 139.56 & 20280 & 45095 \\
\cite{Baran89} (w/ OPF of~\cite{Tse12})  & (10,11), (27,28), (30,31), (7,30), (8,14) & 146.83 &  12.1 &  25.5 \\
\cite{Shirmohammadi89} (w/ OPF of~\cite{Tse12}) & (6,7), (9,10), (13,14), (31,32), (24,28)  & 140.28 & 140.28 & 672 \\
\cite{Gomes06} (w/ OPF of~\cite{Tse12})  & (6,7), (9,10), (13,14), (31,32), (24,28)  & 140.28 & 72.2 &  673 \\
\cite{Khodr09} & (6,7), (8,9), (13,14), (31,32), (24,28)& 139.56  & -- & -- \\
Actual network~\cite{Baran89} & (7,20), (8,14), (11,21), (17,32), (24,28)  & 202.68 & --& --
\end{tabular}}
\end{center}
\label{tab:case2_comparison}
\end{table*}%

\textcolor{black}{Comparisons with an exhaustive search strategy, and also with the heuristic schemes in~\cite{Shirmohammadi89} and~\cite{Gomes06} is provided in Table~\ref{tab:case1_comparison} in terms of computational time and obtained cost. Three points are worth mentioning at this point: \emph{i)} the schemes of~\cite{Shirmohammadi89} and~\cite{Gomes06} were originally developed for balanced systems, and here they have been modified here to account for three-phase lines~\cite{Zidan11} [cf.~\eqref{line_currents}]; specifically, at each iteration, indexes $\sum_{\phi \in \cP_{mn}} |I_{mn}^\phi|^2$ are compared in order to select the switch to open. \emph{ii)} The objective in~\cite{Shirmohammadi89} and~\cite{Gomes06} is to find a tree configuration; thus, the parameter $\lambda$ is set to $200$ in order to obtain a radial network. \emph{iii)} Since~\cite{Shirmohammadi89} and~\cite{Gomes06}  require solving an OPF problem per tested switch status,  the SDP-based reformulation of~\cite{Dallanese-TSG13} is used to solve the OPF optimally. Once the optimal topology is found, OPF is employed to fine-tune voltages, currents, and powers supplied by DG units.}

\textcolor{black}{Complexity is quantified by: \emph{a)} the computational time required by the solver \texttt{SeDuMi}; and, \emph{b)} the overall cpu time, given by the sum of the time required by \texttt{CVX} to first pre-process the data, and by \texttt{SeDuMi} to solve the optimization problems. As expected, the computational time required by the proposed scheme is markedly lower than the competing alternatives [cf.~Remark~2]. As for the optimization objective, the proposed method yields a topology with slightly higher  power loss compared to exhaustive search and~\cite{Shirmohammadi89,Gomes06}. Overall, the power loss in the configurations obtained with the considered methods  is lower than that in the original system. }

\subsection{\textcolor{black}{33-node test system in~\cite{Baran89}}}
\label{sec:feeder2}

\textcolor{black}{The proposed method is tested on the 33-node test system in~\cite{Baran89}, which is broadly considered in the literature for comparison purposes. This is a single-phase system with nominal voltage $12.66$ kV,  37 branches (including tie lines), and total substation loads for the base configuration of 5084.26 kW and 2547 kVAr. No DG units are present in this system. Throughout this subsection, node numbering corresponds to the one in~\cite{Baran89}. Loads are modeled as constant-PQ loads, and transformer losses are neglected (as in~\cite{Baran89, Shirmohammadi89, Gomes06, Khodr09}).}

\textcolor{black}{The sparsity-tuning parameters are set sufficiently high so that a radial topology is obtained. Specifically, two setups are considered: \emph{i)} $\lambda = 2 \times 10^2$ for all lines; and, \emph{ii)} a weighted regularization function is adopted with $\lambda_{mn} = 2 \times 10^3$ for lines $(6,7), (8,9), (9,10), (13,14), (31,32), (7,20), $ $(8,14), (11,21), (17,32), (24,28)$, and $\lambda_{mn} = 2 \times 10^2$ for all the other lines. The second setup represents the case where prior information on the switches that are likely to be opened is available from historical data. The proposed method is compared with~\cite{Baran89, Shirmohammadi89, Gomes06, Khodr09}, as well as with the exhaustive search. Table~\ref{tab:case2_comparison} lists the switches that each scheme suggests to open, the obtained power loss, as well as the required computational time for the consider methods.  It can be seen that the proposed method incurs the lowest complexity. The computational time required by the method in~\cite{Khodr09} is not reported, since it employes a commercial solver that is not publicly available; however, its complexity is expected to be higher than the proposed approach since~\cite{Khodr09} involves the solution of multiple OPF problems (see also Remark~2). When the same $\lambda$ is used for all lines, the proposes method, as well as~\cite{Shirmohammadi89} and~\cite{Gomes06}, outperform~\cite{Baran89} in terms of achieved power loss. A lower power loss is obtained using the method in~\cite{Khodr09}, but at a possibly higher complexity. When a weighted regularization function is used, the proposed method matches the result of~\cite{Khodr09} and that of exhaustive search.}

\subsection{\textcolor{black}{70-node test system in~\cite{Das06}}}
\label{sec:feeder3}

\textcolor{black}{Consider now the $70$-node test system in~\cite{Das06}, which is also used in the literature for comparison purposes. This is a $11$-kV balanced distribution network with two substations, 
four feeders, and 78 branches (including tie lines, with open switches in normal conditions). No DG units are present, and the base topology can be found in~\cite{Das06}. Line parameters as well as load data can be also found in~\cite{Das06}. Specifically, similar to~\cite{Khodr09,Das06}, loads are assumed constant-PQ,  transformer losses are neglected, and the minimum voltage magnitude is set to 0.9 pu.}  

\textcolor{black}{Table~\ref{tab:case3_comparison} summarizes the obtained power losses after reconfiguring the system. It can be seen that the proposed method yields the same topology as~\cite{Shirmohammadi89}  and~\cite{Khodr09}, and outperforms~\cite{Das06}. A comparison between the computational times in Tables~\ref{tab:case2_comparison} and~\ref{tab:case3_comparison} clearly reveals that the proposed method scales well with the network size. In fact, although the number of nodes and lines have doubled, the computational time is approximately the same. This is not the case for~\cite{Shirmohammadi89}, where the computational complexity of OPF solvers grows faster as the network size increases; see also~\cite{PaudyalyISGT} and~\cite{Khodr09} for related claims. }

\begin{table}[t]
\caption{\textcolor{black}{70-node network in~\cite{Das06}: comparison with other methods}}
\begin{center}
\textcolor{black}{\begin{tabular}{|c||c|c|c|}
Method  &  $P_{\mathrm{loss}}$ [kW] & Time for solver [s] &  Overall time [s] \\
\hline
Proposed  & 301.6 & 0.4 & 2.8 \\
\cite{Shirmohammadi89}  & 301.6 & 22.4 & 45.7 \\
\cite{Das06} & 306.9 & -- & -- \\
\cite{Khodr09} & 301.6  & -- & -- \\
Actual network & 341.4  & -- & -- \\
\end{tabular}}
\end{center}
\label{tab:case3_comparison}
\end{table}%

\section{Concluding Remarks}
\label{sec:conclusions}

A DSR problem was considered for three-phase distribution systems featuring DG. Leveraging the notion of group-sparsity, and adopting an approximate linear relation between powers and injected currents, a novel convex DSR formulation was proposed. Being convex, the proposed DSR problem can be solved efficiently even for distribution networks of large size. The ability of the proposed scheme to select the topologies that minimize the overall active power loss was demonstrated via numerical tests, and it was also justified analytically.    

\appendix

\textcolor{black}{
\emph{Proof of Proposition}~\ref{prop:current_switch}. Consider first the case of single-phase distribution lines, where~\eqref{Optimal_current} boils down to
\begin{align}
& \bxi_{mn}^{opt} = \arg \min_{\bxi_{mn} } \, \frac{1}{2} \left(\Re\{Z_{mn}^\phi\}  + \rho_{mn}^{\phi, opt} \right)  \bxi_{mn}^\cT \bxi_{mn}  \nonumber \\
& \hspace{3.5cm} - \bmu_{mn}^\cT  \bxi_{mn}  + \lambda \|\bxi_{mn} \|_2  \, .\hspace{-0.1cm}  \label{Optimal_current_proof}  
\end{align}
The solver of \eqref{Optimal_current_proof} takes the form $\bxi_{mn}  = z \, \bmu_{mn}$ for some scalar $z \geq 0$. In fact, among all possible $\bxi_{mn} $ with the same $\ell_2$-norm, the Cauchy-Schwarz inequality implies that the maximizer of $\bmu_{mn}^{\cT}  \bxi_{mn} $ is colinear with (and in the same direction of) $\bmu_{mn}$. Thus, 
substituting $\bxi_{mn} = z \bmu_{mn}$ into \eqref{Optimal_current_proof} yields the following problem in the scalar $z$:
\begin{align}
& z^{opt} = \arg \min_{ z \geq 0 } \, \frac{1}{2} \left(\Re\{Z_{mn}^\phi\}  + \rho_{mn}^{\phi, opt} \right) z^2 \|\bmu_{mn}\|_2   \nonumber \\
& \hspace{3.5cm} - z \|\bmu_{mn}\|_2  + \lambda |z|  \hspace{-0.1cm}  \label{Optimal_current_proof2}  
\end{align}
where $\lambda |z|$ can be replaced by $\lambda z$ since $z \geq 0$. The necessary and sufficient condition for $z$
to minimize~\eqref{Optimal_current_proof2}  is~\cite[p.~92]{Ruszczynski06}
\begin{equation}     
\left\{
\begin{array}{ll}
   \|\bmu_{mn}\|_2  \leq \lambda,  & \mathrm{ if }\,\, z^{opt} = 0   \\
   \Re\{Z_{mn}^\phi\}  + \rho_{mn}^{\phi, opt}   - \|\bmu_{mn}\|_2 + \lambda = 0, &  \mathrm{ if }\,\, z^{opt} \neq 0    
\end{array}
\right.
\end{equation}
which is satisfied by 
\begin{equation}     
z^{opt} =  \frac{1}{(\Re\{Z_{mn}^\phi\}  + \rho_{mn}^{\phi, opt})\|\bmu_{mn}\|_2)} [\|\bmu_{mn}\|_2 - \lambda]_+ \, .
\end{equation}
}
\textcolor{black}{
Relations~\eqref{Optimal_current_switch2}--\eqref{Optimal_current_eta} can be proved along the lines of~\cite{Wiesel11}.
Specifically,~\eqref{Optimal_current} is first equivalently reformulated as the following quadratic program with a second-order conic constraint
\begin{align}
\hspace{-.2cm} &   \min_{\bxi_{mn}, t } \, \frac{1}{2} \bxi_{mn}^\cT \tbmZ_{mn} \bxi_{mn} + \lambda \|\bxi_{mn} \|_2 - \bmu_{mn}^\cT  \bxi_{mn} + t \nonumber \\
& \mathrm{subject~to} \quad 
\left[ 
\begin{array}{l}
- \lambda   \bxi_{mn} \\
- t    
\end{array}
\right] \preceq 0 \, .
\label{Optimal_current_proof2}  
\end{align}
Next, derive the (concave) dual problem of~\eqref{Optimal_current_proof2}, which amounts to
\begin{align}
\hspace{-.2cm} &   \max_{\bchi } \, - \frac{1}{2} (\bmu_{mn} + \lambda \bchi)^\cT \tbmZ_{mn}^{-1}  (\bmu_{mn} + \lambda \bchi)  \nonumber \\
& \mathrm{subject~to} \quad \|\bchi \|_2^2 \leq 1
\label{Optimal_current_proof3}  
\end{align}
where the constraint $\bchi \in \mathrm{range}(\tbmZ_{mn})$ is left implicit, and $\bchi$ is the multiplier associated with the conic constraint in~\eqref{Optimal_current_proof2}. Consider then the Lagrange dual of~\eqref{Optimal_current_proof3}, namely 
\begin{align}
&  \min_{\eta \geq 0} \max_{\bchi } \, - \frac{1}{2} (\bmu_{mn} + \lambda \bchi)^\cT \tbmZ_{mn}^{-1}  (\bmu_{mn} + \lambda \bchi)  -  \eta \|\bchi \|_2^2 + \eta \, .
\label{Optimal_current_proof4}  
\end{align}
Recalling that $\tbmZ_{mn}$ is invertible, it turns out that the optimal solution of~\eqref{Optimal_current_proof4}  is given by
\begin{align}
\bchi = - \frac{\lambda}{2} \left(\eta \tbmZ_{mn} + \frac{\lambda^2}{2} \bI \right)^{-1} \bmu_{mn}
\label{Optimal_current_proof5}  
\end{align}
with $\eta$ as in~\eqref{Optimal_current_eta}. Notice that the eigenvalues of $(\eta^{opt} \tbmZ_{mn} + \frac{\lambda^2}{2} \bI)^{-1}$, denoted by $\{\theta_i\}$ for brevity, satisfy the inequality $0 < \theta_i \leq 2/\lambda^2$ for all $i = 1,\ldots, 2|\cP_{mn}|$. Thus, when $\|\bmu_{mn}\|_2 \leq \lambda$, the term $(1/2) \bmu_{mn}^\cT \left( \eta \tbmZ_{mn}  + \frac{\lambda^2}{2} \bI_{2|\cP_{mn}|} \right)^{-1} \bmu_{mn}$ in~\eqref{Optimal_current_eta} is negative, thus implying that $\eta^{opt} = 0$. 
}

\bibliographystyle{IEEEtran}
\bibliography{biblio.bib}

\end{document}